\documentclass[12pt]{amsart}
\newcommand{\pof}{\noindent{\em Proof: }}
\usepackage{color}
\newcommand{\s}[1]{\mathcal{#1}}

\newcommand{\C}{\tiny\operatorname{C}}

\newcommand{\M}{{\mathbb M}}

\newtheorem{Thm}{Theorem}[section]
\newtheorem{Def}[Thm]{Definition}
\newtheorem{Rem}[Thm]{Remark}
\newtheorem{Lem}[Thm]{Lemma}

\newtheorem{Prop}[Thm]{Proposition}

\numberwithin{equation}{section}

\fontsize{.5cm}{.5cm}\selectfont\sf

\medskip

\begin{document}

\title[Quantum general linear supergroup]
{Dual canonical bases for the quantum general linear supergroup}
\thanks{Partially supported by the Australian Research Council and Chinese
National Natural Science Foundation project number:10471070}
\author{Hechun Zhang}
\address{Department of Mathematical Sciences, Tsinghua University, Beijing,
China} \email{hzhang@math.tsinghua.edu.cn}
\author{R. B. Zhang}
\address{School of Mathematics and Statistics, University of Sydney,
Australia} \email{rzhang@maths.usyd.edu.au}

\begin{abstract}
Dual canonical bases of the quantum general linear supergroup are
constructed which are invariant under the multiplication of the
quantum Berezinian. By setting the quantum Berezinian to identity,
we obtain dual canonical bases of the quantum special linear
supergroup ${\s O}_q(SL_{m\mid n})$. We apply the canonical bases
to study invariant subalgebras of the quantum supergroups under
left and right translations. In the case $n=1$, it is shown that
each invariant subalgebra is spanned by a part of the dual
canonical bases. This in turn leads to dual canonical bases for
any Kac module constructed by using an analogue of Borel-Weil
theorem.
\end{abstract}
\maketitle

\section{Introduction}
Crystal bases and canonical bases were introduced by Kashiwara
\cite{k1,k2} and Lusztig \cite{lu1,lu2} in the context of
quantized universal enveloping algebras of symmetrizable Kac-Moody
algebras (including finite dimensional simple Lie algebras) and
the associated quantized function algebras in the early 1990s.
Since then their theories have been extensively developed, leading
to many new developments in representation theory. A natural
problem is to study similar bases for the quantized enveloping
superalgebras and the associated quantized function algebras
\cite{zhangrb}, the cousins of the quantized universal enveloping
algebras. Musson and Zou \cite{mz} constructed a crystal basis for
each finite dimensional irreducible module over the the quantized
enveloping superalgebra of $osp(1|2n)$. In \cite{bkk}, the crystal
basis for any highest weight module in a subcategory ${\s
O}_{int}$ of finite dimensional modules over $U_q(gl_{m\mid n})$
is constructed, also see \cite{zou}.  Zou \cite{zou2} also
constructed crystal bases for highest weight modules over the
quantized universal superalgebra of the simple Lie superalgebra
$D(2, 1; \alpha)$ (in Kac's notation).

In contrast, little seems to be known about canonical bases. The
purpose of this paper is to study dual canonical bases for the
quantum general linear supergroup $GL_q{(m\mid n})$ and the
quantum special linear supergroup $SL_q({m\mid n})$. There are
various ways to approach these quantum supergroups. Manin
\cite{manin1, manin2} introduced the coordinate algebra of a
quantum supermatrix and its inverse supermatrix. This superalgebra
has a Hopf superalgebraic structure, thus is regarded as a version
of the quantum general linear supergroup (which depends on more
than one deformation parameters). The paper \cite{zhangrb} studied
the Hopf subalgebra of the finite dual of the quantized universal
enveloping superalgebra $U_q({\mathfrak{gl}}_{m|n})$ generated by
the matrix elements of the natural representation and its dual. It
was shown that this Hopf subalgebra separate points of
$U_q({\mathfrak{gl}}_{m|n})$ in the sense that if $x, y\in
U_q({\mathfrak{gl}}_{m|n})$ are not equal, then there exists $f$
in the Hopf subalgebra such that its evaluations on $x$ and $y$
are different. This Hopf subalgebra was taken as the definition of
another version of the quantum general linear supergroup. It is
not difficult to show that by specializing the parameters of the
Hopf superalgebra of \cite{manin1, manin2}, one obtains the
quantum general linear supergroup of \cite{zhangrb}. A third
version was defined in \cite{f} by localizing the coordinate
algebra of a supermatrix (without the inverse) at the quantum
determinants of two sub-matrices. As the quantum determinants of
these sub-matrices do not commute with each other, properties of
this localization is not immediately transparent. One of the
results of this paper is to show that the third version of the
quantum general linear supergroup is equivalent to the first two.

Any basis for the quantum general linear supergroup will be very
useful for studying its structure.  In \cite{hz}, a basis of the
letter-place algebra (which is a generalization of coordinate
algebra of a quantum supermatrix) was constructed by introducing
quantum minors. In principle one may try to extend this basis to a
basis for the coordinate algebra of the quantum supermatrix
together with its inverse quantum supermatrix, thus to obtain a
basis for the quantum general linear supergroup. However, as far
as we are aware, this was not achieved before, presumably because
of technical difficulties. In fact, it seems that no basis of any
kind is known for the quantum general linear supergroup. In this
paper we shall construct bases for the quantum general linear
supergroup and the associated quantum special linear supergroup.
The results are given in Theorem \ref{basis}.

A notable feature of the bases for the quantum general linear
supergroup is that they are invariant under the multiplication of
the quantum Berezinian. By setting the quantum Berezinian to $1$
we get bases for the quantum special linear supergroup. Also, the
basis elements consist of ${\mathbb Z}[q]$ combinations of certain
monomials, and are invariant under some bar-involution on the
quantum supergroups. Therefore, we may regard the bases
constructed in Theorem \ref{basis} as some dual canonical bases
for the quantum supergroups.

The algebra of functions on the quantum general linear supergroup
admits two actions of the quantized universal enveloping
superalgebra $U_q({\mathfrak{gl}}_{m|n})$. For any bi-subalgebra
$U_S$ of $U_q({\mathfrak{gl}}_{m|n})$, the subspace of invariants
under the left or right translation with respect to $U_S$ forms a
subalgebra, which may be regarded as the algebra of functions on
some quantum homogeneous superspace \cite{zhangrbj} in the general
spirit of non-commutative geometry.  We apply the dual canonical
bases to study such invariant subalgebras. In the case $n=1$, we
show that any subalgebra of invariants is spanned by a part of the
dual canonical bases. This in turn leads to dual canonical bases
for any Kac module constructed using a Borel-Weil type of
construction \cite{zhangrbj}.

The paper is organized as follows. In Section \ref{QUEA-QG}, we
collect some results on the quantized enveloping algebra
$U_q({\mathfrak{gl}}_{m|n})$, and the version of the quantum
general linear supergroup defined in \cite{zhangrb}. The material
of this section will be used throughout the remainder of the
paper. In Section \ref{CA}, the coordinate algebra ${\s
O}_q(M_{m\mid n})$ is presented by exhibiting its generators and
defining relations following \cite{manin1, manin2}; quantum minors
\cite{hz} in the context of ${\s O}_q(M_{m\mid n})$ are discussed;
and the equivalence of the various versions of the quantum general
linear supergroup is proven. In section \ref{BASES}, we present
the construction of the dual canonical bases for ${\s
O}_q(GL_{m\mid n})$ and ${\s O}_q(SL_{m\mid n})$,  and finally in
Section \ref{INV}, we use the dual canonical bases to study
subalgebras of the quantum general linear supergroup which are
invariant under right (or left) translations of any bi-subalgebra
$U_S$ of $U_q({\mathfrak{gl}}_{m|n})$.

\section{Quantum General Linear Supergroup}\label{QUEA-QG}
\subsection{The quantized enveloping algebra
$U_q({\mathfrak{gl}}_{m|n})$}\label{QUEA}

Throughout the paper, we will denote by $\mathfrak g$ the complex
Lie superalgebra ${\mathfrak{gl}}_{m|n}$, and by $U_q({\mathfrak
g})$ the quantized enveloping superalgebra of $\mathfrak g$. Let
${\bf I}=\{1, 2, ..., m+n\}$ and ${\bf I}'={\bf
I}\backslash\{m+n\}$. The quantized enveloping superalgebra
$U_q({\mathfrak g})$ is a ${\mathbb Z}_2$-graded associative
algebra (i.e., associative superalgebra) over ${\mathbb C}(q)$,
$q$ being an indeterminate, generated by $\{K_a, \ K_a^{-1}, \
a\in {\bf I}; \ E_{b\ {b+1}},$ $ \ E_{b+1, b}, \ b\in {\bf I}'\}$,
subject to the relations \cite{zhangrb}
\begin{eqnarray*}
K_a K_a^{-1}=1,
& & K_a^{\pm  1} K_b^{\pm 1} = K_b^{\pm 1}  K_a^{\pm 1}, \nonumber \\
K_a E_{b\ b\pm 1} K_a^{-1} &=&
q_a^{2\delta_{a b} -2\delta_{a\  b\pm 1}} E_{b\  b\pm 1}, \nonumber \\
{}[E_{a\,  a+1},\,E_{b+1\,  b}\}& =& \delta_{a b}
(K_a K_{a+1}^{-1} - K_a^{-1} K_{a+1})/(q_a^2 - q_a^{-2}),\nonumber \\
(E_{m\, m+1})^2 &=& (E_{m+1\, m})^2 = 0, \nonumber \\
E_{a\,  a+1} E_{b\,  b+1} &=& E_{b\,  b+1} E_{a\,  a+1},\nonumber \\
E_{a+1\, a} E_{b+1\, b} &=&E_{b+1\, b} E_{a+1\, a}, \ \ \
\vert a - b\vert \ge 2, \nonumber \\
{\s S}^{(+)}_{a , a\pm 1}&=&{\s S}^{(-)}_{a , a\pm 1}=0,
\ \ , a\ne m,\nonumber \\
E_{m-1\, m+2} E_{m\, m+1} &+& E_{m\, m+1} E_{m-1\, m+2}=0, \\
E_{m+2\, m-1}  E_{m+1\, m} &+& E_{m+1\, m} E_{m+2\, m-1} = 0,
\end{eqnarray*}
where
\begin{eqnarray*}
{\s S}^{(+)}_{a , a\pm 1}&=&
(E_{a\, a+1})^2  E_{a\pm 1\, a+1\pm 1} - (q^2 +
q^{-2}) E_{a\, a+1} \ E_{a\pm 1\, a+1\pm 1} \ E_{a\, a+1}\\
& +& E_{a\pm 1\, a+1\pm1 }\ (E_{a\, a+1})^2,    \\
{\s S}^{(-)}_{a , a\pm 1}&=&
(E_{a+1\, a})^2\,E_{a+1\pm 1\, a\pm 1} - (q^2 +
q^{-2}) E_{a+1\, a}\ E_{a+1\pm 1\, a\pm 1} \ E_{a+1\, a}\\
&+& E_{a+1\pm 1\, a\pm 1}\ (E_{a+1\, a})^2,
\end{eqnarray*}
and $E_{m-1\, m+2}$ and $E_{m+2\, m-1}$ are the $a=m-1$, $b=m+1$,
cases of the following elements
\begin{eqnarray*}
E_{a\, b} &=& E_{a\, c} E_{c\, b} - q_c^{-2} E_{c\, b} E_{a\, c}, \\
E_{b\, a} &=& E_{b\, c} E_{c\, a} -    q_c^2   E_{c\, a} E_{b\,
c}, \ \  a<c<b.
\end{eqnarray*}
Let $[a]=\left\{\begin{array}{l l}
               0,  & \mbox{if} , a\le m, \\
               1,  & \mbox{if} , a>m.
              \end{array}\right. $
Then $q_a=q^{(-1)^{[a]}}$, and \[[E_{a, a+1}, E_{b+1, b}\}= E_{a,
a+1} E_{b+1, b}-(-1)^{[a]+[a+1]} E_{b+1, b} E_{a, a+1}.\] The
${\mathbb Z}_2$ grading of the algebra is specified such that the
elements $K_a^{\pm 1}$, $\forall a\in {\bf I}$, and $E_{b\, b+1}$,
$E_{b+1\,  b}$,  $b\ne m$, are even, while $E_{m\, m+1}$ and
$E_{m+1\, m}$ are odd.  We shall denote by $n^+$ the subalgebra of
$U_q({\mathfrak g})$ generated by $E_{a,a+1}$ for all $a\le m+n-1$
and by $n^-$ the subalgebra generated by $E_{a+1, a}$ for all
$a\le m+n-1$.

It is well known that $U_q({\mathfrak g})$ has the structure of a
${\mathbb Z}_2$ graded Hopf algebra (i.e., Hopf superalgebra),
with a co-multiplication
\begin{eqnarray*}
\Delta(E_{a\, a+1}) &=& E_{a\,  a+1} \otimes
K_a K_{a+1}^{-1} + 1 \otimes E_{a\, a+1}, \\
\Delta(E_{a+1\, a}) &=& E_{a+1\, a }\otimes 1 + K_a^{-1} K_{a+1}
\otimes E_{a+1\, a}, \\
\Delta(K_a^{\pm 1}) &=&K_a^{\pm 1}\otimes K_a^{\pm 1},
\end{eqnarray*}
co-unit
\begin{eqnarray*}
\epsilon(E_{a\, a+1})&=&E_{a+1\, a}=0, \ \ \forall a\in{\bf I}', \\
\epsilon(K_b^{\pm 1})&=&1,  \ \ \ \forall b\in{\bf I},
\end{eqnarray*}
and antipode
\begin{eqnarray*}
S(E_{a\, a+1}) &=& - E_{a\, a+1} K_a^{-1} K_{a+1}, \\
S(E_{a+1\, a}) &=& - K_a K_{a+1}^{-1}E_{a+1\, a}, \\
S(K_a^{\pm 1}) &=&K_a^{\mp 1}.
\end{eqnarray*}
Sometimes, we also use $E_i$ and $F_i$ to denote $E_{i,i+1}$ and
$E_{i+1,i}$ respectively.

Let $\{\epsilon_a | a\in{\bf I}\}$ be the basis of a vector space
with a bilinear for $ (\epsilon_a,\ \epsilon_b
)=(-1)^{[a]}\delta_{a b}$. The roots of the classical Lie
superalgebra ${\mathfrak{gl}}_{m|n}$ can be expressed as $
\epsilon_a - \epsilon_b,$  $a\ne b, \  a, \, b\in{\bf I}.$  It is
known \cite{zhangjmp} that every finite dimensional irreducible
$U_q({\mathfrak g})$ module is of highest weight type and is
uniquely characterized by a highest weight. For $\lambda
=\sum_{a}\lambda_a \epsilon_a$, $\lambda_a\in{\mathbb Z}$, we
shall use the notation $L(\lambda)$ to denote the irreducible
$U_q({\mathfrak g})$ module with a unique (up to scalar multiples)
vector $v_\lambda\ne 0$ such that
\begin{eqnarray*} E_{a, a+1} v_\lambda&=& 0, \ \  a\in{\bf I}',  \\
       K_b v_\lambda&=& q_b^{2\lambda_{b}} v_\lambda,  \ \ \ b\in{\bf I}.
   \end{eqnarray*}
We shall refer to $\lambda$ as the highest weight of $L(\lambda)$.
Then $L(\lambda)$ is finite dimensional if and only if $\lambda$
satisfies $\lambda_a - \lambda_{a+1}\in{\mathbb Z}_+$, $a\ne m$.
In that case, $L(\lambda)$ has the same weight space decomposition
as that of the corresponding irreducible ${\mathfrak{gl}}_{m|n}$
module with the same highest weight $\lambda$.

The natural $U_q({\mathfrak g})$-module $\mathbb E$ has the
standard basis $\{v_a | a\in{\bf I}\}$, such that
\begin{eqnarray*}
K_a v_b = q_a^{2\delta_{a b}} v_b, & & E_{a , a\pm 1} v_b =
\delta_{b, a\pm 1} v_a.
\end{eqnarray*}
The $U_q({\mathfrak g})$ modules ${\mathbb E}^{\otimes k}$,
$k\in{\mathbb Z}_+$ ( ${\mathbb E}^0= {\mathbb C}$ ) were shown to
be completely reducible \cite{zhangrb}. The irreducible summands
of these $U_q({\mathfrak g})$-modules, referred to as irreducible
contravariant tensor modules, can be characterized in the
following way. Let ${\mathbb Z}_+$ be the set of nonnegative
integers. Define a subset $\s P$ of ${\mathbb Z}_+^{\otimes(m+n)}$
by
\begin{eqnarray*}
{\s P}&=&\{ p=(p_1, p_2, ..., p_{m+n})\in {\mathbb Z}_+^{m+n}
\mid\ p_{m+1}\le n, \  p_a\ge p_{a+1},\   a\in{\bf I}'\}.
\end{eqnarray*}
We associate with each $p\in{\s P}$  a
${\mathfrak{gl}}_{m|n}$-weight
 defined by
\begin{eqnarray*}
\lambda^{(p)}=\sum_{i=1}^{m} p_i \epsilon_i + \sum_{\nu=1}^n
\sum_{\mu=1}^{p_{m+\nu}} \epsilon_{m+\mu},
\end{eqnarray*}
and let
\begin{eqnarray}
\Lambda^{(1)}&=&\{ \lambda^{(p)} \mid  p\in {\s P} \}.
\end{eqnarray}
{}From results of \cite{Gould, I} we know that an irreducible
$U_q({\mathfrak g})$-module is a contravariant tensor if and only
if its highest weight belongs to $\Lambda^{(1)}$.

Let $L(\lambda)$  be  an irreducible contravariant tensor
$U_q({\mathfrak g})$ module with highest weight
$\lambda\in\Lambda^{(1)}$. We define $\bar{\lambda}$ to be its
lowest weight, and set $\lambda^\dagger=-\bar{\lambda}$. An
explicit formula for $\lambda^\dagger$ was given in \cite{Gould}
(section III. B.), where a more compact characterization was also
given for $\Lambda^{(1)}$ and also the set
\begin{eqnarray}
\Lambda^{(2)}&:=& \{\lambda^\dagger\ | \ \lambda\in\Lambda^{(1)}\}.
\end{eqnarray}
We refer to that paper for details. Now  the dual module
$L(\lambda)^\dagger$ of $L(\lambda)$, which we will  call a
covariant tensor module, has highest weight $\lambda^\dagger$. The
most important example is the dual module ${\mathbb
E}^\dagger=L(-\epsilon_{m+n})$ of ${\mathbb E}$.

The situation with tensor powers ${\mathbb E}^{\otimes k}$   and
$({\mathbb E}^\dagger)^{\otimes k}$ can be summarized into the
following proposition \cite{zhangrb}.
\begin{Prop}\label{tensor}
\begin{enumerate}
\item Each $U_q({\mathfrak g})$-module ${\mathbb E}^{\otimes k}$
    (resp. $({\mathbb E}^\dagger)^{\otimes k}$),
   $k\in {\mathbb Z}_+$, can be decomposed into a direct sum of irreducible
   modules with highest weights belonging to $\Lambda^{(1)}$
   (resp. $\Lambda^{(2)}$).
\item Every irreducible $U_q({\mathfrak g})$-module with highest weight belonging
    to $\Lambda^{(1)}$ (resp. $\Lambda^{(2)}$) is a direct summand of
    some tensor powers of ${\mathbb E}$ (resp. ${\mathbb E}^\dagger$).
\end{enumerate}
\end{Prop}

\subsection{The quantum general linear supergroup}\label{QG}

Let $(U_q({\mathfrak g}))^0$ be the finite dual of $U_q({\mathfrak
g})$, which, by standard Hopf algebra theory,  is a ${\mathbb
Z}_2$-graded Hopf algebra with structure dualizing that of
$U_q({\mathfrak g})$. Let us denote by $\pi$ the representation of
$U_q({\mathfrak g})$ on $\mathbb E$ relative to the standard basis
$\{ v_a\mid a\in{\bf I}\}$:
\begin{eqnarray*}x v_a&=& \sum_{b} \pi(x)_{b\, a} v_b, \ \ \
x\in U_q({\mathfrak g}),
\end{eqnarray*}
then we have the elements $t_{a\, b}\in (U_q({\mathfrak g}))^0$,
$a, b\in{\bf I}$, defined by
\begin{eqnarray*}t_{a\, b} ( x ) &=& \pi ( x )_{a\, b}, \ \ \ \forall x\in
U_q({\mathfrak g}). \end{eqnarray*} Note that $t_{a\, b}$ is even
if $[a]+[b]\equiv 0 ( mod\, 2 )$, and odd otherwise.

Consider the subalgebra $G_q^{\pi}$ of $(U_q({\mathfrak g}))^0$
generated by $t_{a\, b}$, $a, b\in{\bf I}$. The multiplication
which $G_q^{\pi}$ inherits from $(U_q({\mathfrak g}))^0$ is given
by
\begin{eqnarray*}\langle t\ t',\ x\rangle &=& \sum_{(x)} \langle t\otimes t',
\ x_{(1)}\otimes x_{(2)} \rangle\nonumber\\
&=& \sum_{(x)} (-1)^{[t'][x_{(1)}]} \langle t,  x_{(1)} \rangle
\langle t',\ x_{(2)} \rangle, \ \ \ \ \forall  t, t'\in G_q^{\pi},
\ x\in U_q({\mathfrak g}).
\end{eqnarray*}
To better understand the algebraic structure of $G_q^\pi$, recall
that the Drinfeld version of $U_q({\mathfrak g})$ admits a
universal $R$ matrix, which in particular satisfies
\begin{eqnarray*} R \Delta( x ) &=&
\Delta'( x ) R,  \ \ \ \forall x\in U_q({\mathfrak g}).
\end{eqnarray*}
Applying $\pi\otimes\pi$ to both sides of the equation yields
\begin{eqnarray}R^{\pi\,
\pi} (\pi\otimes\pi)\Delta( x ) &=& (\pi\otimes\pi)\Delta'( x )
R^{\pi\, \pi}, \label{R} \end{eqnarray} where $R^{\pi\, \pi}
:=(\pi\otimes\pi)R$ is given by
\begin{eqnarray*}
R^{\pi\, \pi} &=&q^{2\sum_{a\in{\bf I}} e_{a\, a}\otimes e_{a\, a}
(-1)^{[a]} }
 + (q^2-q^{-2})\sum_{a<b} e_{a\, b}\otimes e_{b\, a}
(-1)^{[b]} .
\end{eqnarray*}

As we work with the Jimbo version of $U_q({\mathfrak g})$ in this
paper, it is problematic to talk about a universal $R$ matrix.
However, it is important to note that equation (\ref{R}) makes
perfect sense in the present setting.

We can re-interpret equation (\ref{R}) in terms of $t_{a b}$ in
the following way. Set $t=\sum_{a, b} e_{a\, b}\otimes t_{a\, b}$.
Then
\begin{eqnarray}R^{\pi\, \pi}_{1 2}\, t_1\, t_2 &=& t_2\, t_1\, R^{\pi\,
\pi}_{1 2}.  \label{Faddeev} \end{eqnarray} The co-multiplication
$\Delta$ of $G_q^{\pi}$ is also defined in the standard way by
\begin{eqnarray*}
\langle \Delta( t_{a\, b} ), \ x\otimes y \rangle&=& \langle
t_{a\, b} , \ x y\rangle = \pi (x y )_{a\, b}, \ \ \ \ \forall x,
y\in U_q({\mathfrak g}).
\end{eqnarray*}
We have
\begin{eqnarray}\Delta( t_{a\, b} )&=&\sum_{c\in{\bf I}}
(-1)^{([a]+[c])([c]+[b])} t_{a\, c}\otimes t_{c\, b}.
\end{eqnarray}
$G_q^{\pi}$ also has the unit $\epsilon$, and the co - unit
$1_{U_q({\mathfrak g})}$. Therefore, $G_q^{\pi}$ has the
structures of a ${\mathbb Z}_2$-graded bi-algebra.

Let $\pi^{(\lambda)}$ be an arbitrary irreducible contravariant
tensor representation of $U_q({\mathfrak g})$. Define the elements
$t^{(\lambda)}_{i\, j}$, $i, \, j= 1, 2, ..., dim_{\mathbb
C}\pi^{(\lambda)}$, of $(U_q({\mathfrak g}))^0$ by
\begin{eqnarray*}
t^{(\lambda)}_{i\, j}( x )&=& \pi^{(\lambda)}( x )_{i\, j}, \ \ \
\forall x\in U_q({\mathfrak g}).
\end{eqnarray*}
These will be called the matrix elements of the irreducible
representation $\pi^{(\lambda)}$. It is an immediate consequence
of Proposition \ref{tensor} that for every
$\lambda\in\Lambda^{(1)}$, the elements $t^{(\lambda)}_{i\, j}\in
G_q^\pi$, for all $i,\ j$, and every $f\in G_q^\pi$ can be
expressed as a linear sum of such elements. We have the following
result \cite{zhangrb}.
\begin{Prop}\label{peter-weyl}
As a vector space,
\begin{eqnarray*} G_q^\pi=\bigoplus_{\lambda\in\Lambda^{(1)}}
T^{(\lambda)}, \text{ where } T^{(\lambda)}=\bigoplus_{i,
j=1}^{dim\pi^{(\lambda)}} {\mathbb C}(q) t^{(\lambda)}_{i\, j}.
\end{eqnarray*}
\end{Prop}

Let $\{{\bar v}_a\ | \ a\in{\bf I}\}$ be the basis of ${\mathbb
E}^\dagger$ dual to the standard basis of ${\mathbb E}$, i.e.,
\begin{eqnarray*}{\bar v}_a ( v_b )&=\delta_{a \, b}. \end{eqnarray*}
Denote by $\bar{\pi}$ the covariant vector irreducible
representation relative to this basis. Let ${\bar t}_{a\, b}$, $a,
b\in{\bf I}$, be the elements of $(U_q({\mathfrak g}))^0$ such
that
\begin{eqnarray*}{\bar t}_{a\, b}( x )&=&{\bar\pi}( x )_{a\, b}, \
\ \ \forall x\in U_q({\mathfrak g}). \end{eqnarray*}Note that
${\bar t}_{a\, b}$ is even if $[a]+[b]\equiv 0 ( mod\,  2 )$, and
odd otherwise.  These elements generate a ${\mathbb Z}_2$-graded
bi-subalgebra $G_q^{\bar\pi}$ of $(U_q({\mathfrak g}))^0$. We want
to point out that the ${\bar t}_{a\, b}$ obey the relation
\begin{eqnarray}R^{{\bar\pi}\, {\bar\pi}}_{1 2}\ {\bar
t}_1\ {\bar t}_2&=& {\bar t}_2\ {\bar t}_1 R^{{\bar\pi}\,
{\bar\pi}}_{1 2}, \label{YB} \end{eqnarray} where ${\bar
t}=\sum_{a, b} e_{a\, b}\otimes {\bar t}_{b\, a}$ and
$R^{{\bar\pi}\, {\bar\pi}}= ({\bar\pi}\otimes{\bar\pi})R$ is given
by
\begin{eqnarray*}
R^{{\bar\pi}\, {\bar\pi}}&=& q^{2\sum_{a\in{\bf I}} e_{a\, a}\otimes e_{a\, a}
(-1)^{[a]} }
 + (q^2-q^{-2})\sum_{a>b} e_{a\, b}\otimes e_{b\, a} (-1)^{[b]}.
\end{eqnarray*}
Also, the co-multiplication is given by
\begin{eqnarray*}
\Delta( \overline{t}_{a\, b} )&=&\sum_{c\in{\bf I}}
(-1)^{([a]+[c])([c]+[b])} \overline{t}_{a\, c}\otimes \overline{t}_{c\, b}.
\end{eqnarray*}
Similar to the case of $G_q^{\pi}$, we let $\bar{T}^{(\lambda)}$
be the subspace of $G_q^{\bar{\pi}}$ spanned by the matrix
elements of the irreducible representation with highest weight
$\lambda\in \Lambda^{(2)}$. Then it follows from Proposition 2.1
that
\begin{Prop} As a vector space,
\begin{eqnarray*}
G_q^{\bar\pi}&=&\bigoplus_{\mu\in\Lambda^{(2)}} {\bar T}^{(\mu)}.
\end{eqnarray*}
\end{Prop}
\begin{Def} The algebra $G_q$ of functions on the quantum general linear supergroup
$GL_q({m\mid n})$ is the ${\mathbb Z}_2$-graded subalgebra of
$U_q(\mathfrak g)^0$ generated by $\{t_{ab},\overline{t}_{ab}\mid
a,b\in {\bf I}\}$.\end{Def}

The  $t_{a\, b}$ and $\overline{t}_{a\, b}$, besides obeying the
relations (\ref{Faddeev}) and (\ref{YB}), also  satisfy
\begin{eqnarray}R^{{\bar\pi}\, \pi}_{1 2}\ {\bar t}_1 \ t_2 &=&
t_2 \ {\bar t}_1 \ R^{{\bar\pi}\, \pi}_{1 2}, \label{mix}
\end{eqnarray}
where $R^{{\bar\pi}\, \pi} := ({\bar\pi}\otimes \pi) R$ is given
by
\begin{eqnarray*}
R^{{\bar\pi}\, \pi} &=& q^{-2\sum_{a\in{\bf I}} e_{a\, a}\otimes
e_{a\, a} (-1)^{[a]} } - (q^2-q^{-2})\sum_{a<b} e_{b\, a}\otimes
e_{b\, a} (-1)^{[a]+ [b] + [a][b]}.
\end{eqnarray*}
As both $G_q^\pi$ and $G_q^{\bar\pi}$ are ${\mathbb Z}_2$-graded
bi-algebras, $G_q$ inherits a natural bi-algebra structure. It
also admits an antipode $S: G_q\rightarrow G_q$, which is a linear
anti-automorphism given by
\begin{eqnarray}
S(t_{a\, b})= (-1)^{[a][b]+[a]} {\bar t}_{b\, a}, \quad  S({\bar
t}_{a\, b})= (-1)^{[a][b]+[b]} q^{2( 2\rho,\  \epsilon_a -
\epsilon_b )} t_{b\, a} .
\end{eqnarray}
Therefore, $G_q$ has the structures of a ${\mathbb Z}_2$-graded
Hopf algebra.  It was also shown  in \cite{zhangrb} that $G_q$
seperates points of $U_q(\mathfrak g)$ in the sense that for any
$x,y\in U_q(\mathfrak g)$ such that $x\ne y$, then there exists
$f\in G_q$ such that $\langle f,\ x\rangle \ne \langle f,\
y\rangle$.

\medskip

Denote by $-$ the involution of the base field ${\mathbb Q}(q)$
which is given as follows:
\begin{eqnarray*}-: {\mathbb Q}(q)\longrightarrow  {\mathbb
Q}(q),\quad  q\mapsto q^{-1}.\end{eqnarray*} There is an
automorphism of the quantized enveloping algebra $U_q(\mathfrak
g)$ denoted by $\bar{ }$ which is given by:
\[\bar{E_i}=E_i,\quad \bar{F_i}=F_i,
\quad \bar{K_j}=K_j^{-1},\quad \bar{q}=q^{-1}\text{ for all }i,j.\]
We define a  linear map $\dagger$ on $G_q$ as follows:
\[<f^\dagger,x>= <f, \bar{x}>^-,\]
for any $f\in G_q$ and $x\in U_q(\mathfrak g)$.

The Hopf superalgebra $G_q$ admits a natural bi-module structure
over $U_q(\mathfrak g)$, with the left action given below:
\[E_i t_{kl}=\delta_{i, k-1} t_{il},\quad
F_i t_{kl} =\delta_{i, k} t_{i+1,l},\quad
K_it_{kl}=q^{2\delta_{ik}}.\]
The right action can be written down similarly (see section 6 for
more details). Hence, we can talk about the left weight and right
weight with respect to the action of the Cartan subalgebra. An
element is called homogeneous if it is both left weight vector and
right weight vector. For  a homogenous element $f$, denote by
$w_l(f)$ the left weight of $f$ and by $w_r(f)$ the right weight
of $f$.

\begin{Prop}\label{dagger}
For any homogeneous elements $f, g\in G_q$ with weights $w_l(f),
w_r(f)$ and $w_l(g), w_r(g)$ respectively,
\[(fg)^\dagger=(-1)^{[f][g]}q^{2(w_l(f),
w_l(g))- 2(w_r(f), w_r(g))}g^\dagger f^\dagger.\]
\end{Prop}

\pof Denote by $\Theta$ the quasi $R$-matrix which satisfies the
following identities in the completion of $U_q({\mathfrak
g})\otimes U_q({\mathfrak g})$:
\[\Theta\overline{\Theta}=\overline{\Theta}\Theta=1\otimes 1,\]
\[\Theta(-\otimes -)\Delta(x)=\Delta(\bar{x})\Theta.\]
Let $\Phi$ be the algebra automorphism of $U_q({\mathfrak
g})\otimes U_q({\mathfrak g})$ defined by
\begin{eqnarray*}
E_i\otimes 1\mapsto E_i\otimes k_i,&&
1\otimes E_i\mapsto k_i\otimes E_i,\\
F_i\otimes 1\mapsto E_i\otimes k_i^{-1},
&& 1\otimes F_i \mapsto k_i^{-1}\otimes F_i,\\
k_i\otimes1 \mapsto k_i\otimes 1, && 1\otimes k_i \mapsto 1
\otimes k_i,
\end{eqnarray*} for all $i=1,2,\cdots, m+n-1$. One can check
easily that
\[\Delta(\bar{x})=\Phi\circ (-\otimes -)\circ \Delta^\prime(x),\]
where $\Delta^\prime$ is the opposite co-multiplication.

Now, for any homogeneous elements $f$ and $g$, we have
\begin{eqnarray*}<(fg)^\dagger, x>&=&<fg, \bar{x}>^-
=<f\otimes g, \Delta(\bar{x})>^-\\\nonumber
&=&<f\otimes g, \Phi\circ (-\otimes -)\circ \Delta^\prime(x)>^-\\\nonumber
&=& q^{2(w_r(f), w_r(g))- 2(w_l(f), w_l(g))}
   <f\otimes g, (-\otimes -)\Delta^\prime(x)>^-\\\nonumber
&=&q^{2(w_r(f), w_r(g))- 2(w_l(f), w_l(g))}
  <f^\dagger\otimes g^\dagger, \Delta^\prime(x)>\\\nonumber
&=&(-1)^{[f][g]}q^{2(w_r(f), w_r(g))- 2(w_l(f), w_l(g))}
  <g^\dagger f^\dagger, x>.\end{eqnarray*}
\qed

The linear map $\dagger$ can be furnished into an anti
automorphism in the following way:
\begin{Lem}\label{anti} The mapping $\bar{ }:G_q \longrightarrow
G_q$ defined, for any homogeneous element $f\in G_q$ with weights
$(w_l(f), w_r(f))$,  by
\begin{eqnarray*}f\mapsto
q^{(w_l(f),w_l(f)-(w_r(f),w_r(f))}f^\dagger, \quad  q\mapsto
q^{-1}
\end{eqnarray*}
is an anti-automorphism of the superalgebra $G_q$.
\end{Lem}

\pof For any homogeneous elements $f, g\in G_q$ with weights
$(w_l(f), w_r(f))$ and $(w_l(g), w_r(g))$ respectively, we have
\begin{eqnarray}
\overline{fg}&=&q^{A(f, g)} (fg)^\dagger=(-1)^{[f][g]}q^{B(f,
g)}g^\dagger f^\dagger, \label{fg-dagger}\end{eqnarray}
where
\begin{align*}
A(f, g)=&(w_l(f)+w_l(g),w_l(f)+w_l(g))-
    (w_r(f)+w_r(g),w_r(f)+w_r(g)), \\
B(f, g)=&A(f,g)+ 2(w_r(f), w_r(g))- 2(w_l(f), w_l(g)).
\end{align*}
The far right hand side of equation (\ref{fg-dagger}) can be
easily shown to be equal to $(-1)^{[f][g]}\bar{g}\bar{f}$, thus
completing the proof.

\qed

\section{Coordinate algebra  of  quantum supermatrix} \label{CA}

Let $X$ be an $(m+n)\times (m+n)$ quantum supermatrix. We shall
always write it in block form
\[X=\begin{pmatrix}A&B\\C&D\end{pmatrix}\] where $A$ and $D$ are
respectively $m\times m$ and $n\times n$ sub-matrices of even
entries, while $B$ and $C$ are respectively $m\times n$ and
$n\times m$ sub-matrices of odd entries. The coordinate algebra
${\s O}_q(M_{m\mid n})$ \cite{manin1} of $X$ is a ${\mathbb
Z}_2$-graded algebra (i.e., superalgebra) generated by the entries
of the quantum supermatrix $X$ satisfying the relation
$$R^{\pi, \pi}X_1X_2=X_2X_1R^{\pi, \pi},$$
where $X_1=X\otimes 1$ and $X_2=1\otimes X$. The defining
relations can be written explicitly as follows:
\begin{eqnarray*}\label{srelation}
x_{ij}x_{ik}&=&(-1)^{([i]+[j])([i]+[k])}
           q^{2(-1)^{[i]}}x_{ik}x_{ij}, j<k,\\\nonumber
x_{ij}x_{kj}&=&(-1)^{([i]+[j])([k]+[j])}
            q^{2(-1)^{[j]}}x_{kj}x_{ij}, i<k,\\\nonumber
x_{ij}x_{kl}&=&(-1)^{([i]+[j])([k]+[l])}x_{kl}x_{ij}, i<k,
j>l,\\\nonumber
x_{ij}x_{kl}&=&(-1)^{([i]+[j])([k]+[l])}x_{kl}x_{ij}+
(-1)^{[k][j]+[k][l]+[j][l]} (q^2-q^{-2})x_{il}x_{kj}, \\
&&i<k, j<l.\end{eqnarray*} Note that the matrix $A$ is a quantum
matrix with deformation parameter $q$ while $D$ is also a quantum
matrix with deformation parameter $q^{-1}$. We shall refer $A$ as
a $q$-matrix and $D$ as a $q^{-1}$-matrix. The superalgebra ${\s
O}_q(M_{m\mid n})$ has a bi-superalgebra structure, and it is
customary to take the following coproduct and counit:
\begin{eqnarray*}
\Delta(x_{ij})=\sum_k x_{ik}\otimes x_{kj}, &\quad&
\epsilon(x_{ij})=\delta_{ij}.
\end{eqnarray*}

Let us recall the definition of quantum minors in the present
context. We shall largely follow \cite{hz}, besides some slight
change in conventions. Denote 
\begin{eqnarray*} (\ , \ ):{\mathbb
Z}^{m+n}\times {\mathbb Z}^{m+n}\longrightarrow {\mathbb
Z},\\\nonumber (\underline{a},
\underline{b})=\sum_{i>j}(a_ib_j-a_jb_i),
\end{eqnarray*}
where $\underline{a}=(a_1,a_2,\cdots, a_{m+n}),
\underline{b}=(b_1, b_2,\cdots, b_{m+n})\in {\mathbb Z}  ^{m+n}$.
Also, denote by $$[\underline{a}]=\sum_i a_i[i]\text{ and }
\{\underline{a}\}=\sum_i a_i([i]+1).$$

The quantum superspace (or rather its coordinate algebra) $A_q$ is
a superalgebra generated by $x_1,x_2,\cdots, x_m, \cdots, x_{m+n}$
with parity assignment $[x_i]=[i]$ and defining relations
\begin{eqnarray*}x_i^2&=&0, \quad \text{if }[i]=\bar{1},\\\nonumber
x_ix_j&=&(-1)^{[i][j]}q^2x_jx_i, \quad \text{for }
i<j.\end{eqnarray*}
We also introduce the superalgebra $A_q^*$
generated by $\xi_1,\xi_2,$ $\cdots,$ $\xi_m,$ $\cdots,$
$\xi_{m+n}$ with parity assignment $[\xi_i]=\bar{1}-[i]$ and
defining relations
\begin{eqnarray*}
\xi_i^2&=&0, \quad \text{ if }[i]=\bar{0},\\ \nonumber
\xi_i\xi_j&=&(-1)^{([i]+\bar{1})([j]+\bar{1})}q^2\xi_j\xi_i, \quad
\text{ for }i>j.
\end{eqnarray*}
For any $\underline{a}=(a_1,a_2,\cdots, a_{m+n})\in {\mathbb
Z}_+^{m+n}$, denote by
$$x^{\underline{a}}=x_1^{a_1}x_2^{a_2}\cdots
x_{m+n}^{a_{m+n}},\quad \xi^{\underline{a}}=\xi_1^{a_1}
\xi_2^{a_2}\cdots \xi_{m+n}^{a_{m+n}},$$ the ordered monomials.
Then $x^{\underline{a}}$, $\overline{a} \in {\mathbb Z}_+^m\times
{\mathbb Z}_2^n$, form a basis of $A_q$, while
$\xi^{\overline{b}}$, $\overline{b}\in {\mathbb Z}_2^m\times
{\mathbb Z}_+^n$,  form a basis of $A_q^*$.

The superalgebras $A_q$ and $A_q^*$ are comodule superalgebras of
${\s O}_q(M_{m\mid n})$ with the coactions given below:
\begin{Thm}(\cite{manin1}) There exist superalgebra morphisms
\begin{eqnarray*}
\delta: A_q\longrightarrow  {\s O}_q(M_{m\mid n})\otimes
A_q, &\quad& x_i\mapsto \sum_j x_{ij}\otimes x_j;\\
\delta^*: A_q^*\longrightarrow  {\s O}_q(M_{m\mid n})\otimes A_q^*
&\quad& \xi_i \mapsto \sum_j x_{ij}\otimes \xi_j,\end{eqnarray*}
which give ${\s O}_q(M_{m\mid n})$-comodule structures to $A_q$
and $A_q^*$.
\end{Thm}

Define elements $\Delta(\underline{a}, \underline{b})$  and
$\Delta(\underline{a}, \underline{b})^*$ of ${\s O}_q(M_{m\mid
n})$ by
\[\delta(x^{\underline{a}})=\sum_{\underline{b}}\Delta(\underline{a},
\underline{b})\otimes x^{\underline{b}},\qquad
\delta^*(\xi^{\underline{a}})=\sum_{\underline{b}}\Delta(\underline{a},
\underline{b})^*\otimes \xi^{\underline{b}},\] which will be
referred to as quantum minors. Note that our definition differs
from that of \cite{hz} by a scalar. See also \cite{nmy, pw} for
discussions of quantum minors.

The following quantum Laplace expansion can be derived directly
from the definition of the quantum minors.
\begin{Prop}
\begin{eqnarray*}
\Delta(\underline{a}+\underline{a}^\prime, \underline{b})&=&
\sum_{\underline{b}=\underline{c}+\underline{c}^\prime}
(-1)^{[\underline{c}]
([\underline{a}^\prime]+[\underline{c}^\prime])}
q^{2(\underline{a}, \underline{a}^\prime)-
2(\underline{c},\underline{c}^\prime)}\Delta(\underline{a},
\underline{c})\Delta( \underline{a}^\prime,
\underline{c}^\prime);\\
\Delta(\underline{a}+\underline{a}^\prime, \underline{b})^*&=&
\sum_{\underline{b}=\underline{c}+\underline{c}^\prime}
(-1)^{\{\underline{c}\}
(\{\underline{a}^\prime\}+\{\underline{c}^\prime\})}
q^{2(\underline{a}, \underline{a}^\prime)-
2(\underline{c},\underline{c}^\prime)}\Delta(\underline{a},
\underline{c})^*\Delta( \underline{a}^\prime,
\underline{c}^\prime)^*.\end{eqnarray*}
\end{Prop}

\pof Consider $\Delta(\underline{a}, \underline{b})$. We have
\begin{eqnarray*}
\delta(x^{\underline{a}+\underline{a}^\prime})&=&q^{2(\underline{a},
\underline{a}^\prime)}
\delta(x^{\underline{a}})\delta(x^{\underline{a}^\prime})\\\nonumber
&=&q^{2(\underline{a}, \underline{a}^\prime)}(\sum_{\underline{c}}
\Delta(\underline{a}, \underline{c})\otimes x^{\underline{c}})
(\sum_{\underline{c}^\prime} \Delta(\underline{a}^\prime,
\underline{c}^\prime)\otimes x^{\underline{c}^\prime})\\\nonumber
&=&\sum_{\underline{b}}
\sum_{\underline{c}+\underline{c}^\prime=\underline{b}}
(-1)^{[\underline{c}]
([\underline{a}^\prime]+[\underline{c}^\prime])}
q^{2(\underline{a}, \underline{a}^\prime)-
2(\underline{c},\underline{c}^\prime)}\Delta(\underline{a},
\underline{c})\Delta( \underline{a}^\prime,
\underline{c}^\prime)\otimes x^{\underline{b}}.
\end{eqnarray*}
The quantum Laplace expansion for $\Delta(\underline{a},
\underline{b})^*$ can be proved similarly.

\qed

For any $r,s\in {\mathbb Z}, r<s$, denote by $[r, s]=\{r,
r+1,\cdots,s\}$. The following quantum minors will play important
roles later. Let
\begin{eqnarray*}
det_q A&:=&\Delta([1,m],[1,m])^*, \\
det_{q^{-1}}D&:=&\Delta([m+1,m+n], [m+1,m+n]).
\end{eqnarray*}
Then
\begin{eqnarray*}
det_q A &=&\sum_{\sigma\in
S_m}(-q^2)^{l(\sigma)}x_{1\sigma(1)}x_{2\sigma(2)}\cdots
x_{m\sigma(m)}, \\
det_{q^{-1}} D &=&\sum_{\tau\in S_n}
(-q^{-2})^{l(\tau)}x_{m+1,m+\tau(1)}x_{m+2,m+\tau(2)}\cdots
x_{m+n,m+\tau(n)}.
\end{eqnarray*}
For $r\le min\{m,n\}$, we have
\begin{eqnarray*}
\Delta([1,r],[m+1, m+r])^* &=& \sum_{\sigma\in
S_r}(-q^2)^{l(\sigma)}x_{1,m+\sigma(1)}x_{2,m+\sigma(2)}\cdots
x_{r,m+\sigma(r)},\\
\Delta([m+1,m+r],[1,r])&=&\sum_{\sigma\in
S_r}(-q^{-2})^{l(\sigma)}x_{m+1,\sigma(1)}x_{m+2,\sigma(2)}\cdots
x_{m+r,\sigma(r)}. \end{eqnarray*} As we shall see later, the
quantum minors $det_q A$ and $\Delta([1,r],[m+1, m+r])^*$ are
annihilated by every $E_i$ under left translation,  while
$det_{q^{-1}}D$ and $\Delta([m+1,m+r],[1,r])$ are annihilated by
every $F_i$ under left translation.  However, the following minors
\begin{eqnarray*} &&\Delta([1,m],[1,m])=\sum_{\sigma\in
S_m}q^{-2l(\sigma)}x_{1\sigma(1)}x_{2\sigma(2)}\cdots
x_{m\sigma(m)}, \\
&&\Delta([m+1,m+n], [m+1,m+n])^*
\\ && =\sum_{\tau\in S_n}
q^{-2l(\tau)}x_{m+1,m+\tau(1)}x_{m+2,m+\tau(2)}\cdots
x_{m+n,m+\tau(n)}\end{eqnarray*}
do not behave well under the left
and right translations.

Let ${\s S}$ be the multiplicative set of products of powers of
$det_qA$ and $det_{q^{-1}}D$. Denote by ${\s O}_q(GL_{m\mid n})$
the localization of ${\s O}_q(M_{m\mid n})$ at ${\s S}$. Then the
inverse matrix $X^{-1}$ of $X$ lies in ${\s O}_q(GL_{m\mid n})$,
as can be shown by an explicit calculation \cite{phung}. We shall
always write $X^{-1}$ in the block form
\[X^{-1}=\begin{pmatrix}\bar{A}&\bar{B}\\\bar{C}&\bar{D}\end{pmatrix},\]
where $\bar{A}$ is an $m\times m$ $q^{-1}$-matrix while $\bar{D}$
is a $n\times n$ $q$-matrix \cite{zhangrb}. In \cite{ls}, the
explicit formula for the quantum Berezinian was given:
\[Ber_q=det_qA det_{q^{-1}}\bar{D},\]
which is also known to be central in ${\s O}_q(GL_{m\mid n})$.

\begin{Rem}
The commutation relations among all of the entries in the matrices
$X$ and $X^{-1}$ were given in \cite{zhangrb} by using $R$
matrices.
\end{Rem}

Let us define the following quantum $q^{-1}$-matrix,
\[ D^\prime:=(y_{\mu,\nu})_{\mu,\nu=m+1}^{m+n}
=D-CA^{-1}B, \] the entries of which all belong to  ${\s
O}_q(GL_{m\mid n})$.  Also note that $D^\prime$ is in fact the
inverse matrix of $\bar{D}$.

\begin{Prop} \label{commun} The following commutation relations hold:
\begin{eqnarray*}
det_qA x_{ij}=x_{ij}det_qA, &\quad& det_qA x_{ \mu,j}=q^2 x_{
\mu,j}det_qA,\\\nonumber
det_qA x_{i,  \nu}=q^2  x_{i, \nu}det_qA, &\quad& det_qA
y_{\mu,\nu}=y_{\mu, \nu}det_qA;\\\nonumber
det_{q^{-1}}D^\prime x_{ij}=x_{ij} det_{q^{-1}}D^\prime,&\quad&
det_{q^{-1}}D^\prime x_{ \mu,j}=q^{2} x_{
\mu,j}det_{q^{-1}}D^\prime, \\\nonumber
det_{q^{-1}}D^\prime x_{i,  \nu}=q^{2} x_{i,
\nu}det_{q^{-1}}D^\prime, &\quad& det_{q^{-1}}D^\prime y_{\mu,
\nu}=y_{\mu, \nu}det_{q^{-1}}D^\prime
\\\nonumber
\text{for } i,j=1,2,\cdots,m,  &\quad& \mu,\nu=m+1,m+2,\cdots,
m+n.\end{eqnarray*} Moreover, $det_qA(det_{q^{-1}}D^\prime)^{-1}$
is central in ${\s O}_q(GL_{m\mid n})$.
\end{Prop}

\pof The inverse matrix of the quantum matrix $A$ is
$$A^{-1}=\left((-q^2)^{j-i}A_{ji}(det_qA)^{-1}\right),$$
where $A_{ji}$ is the determinant of the quantum matrix obtained
from $A$ by deleting the $j$th row and $i$th column. Now the claim
that $det_qA$ commutes with $y_{\mu, \nu}$ is equivalent to the
relation
\[det_qA x_{\mu, \nu}= x_{\mu, \nu}det_qA +(q^2-q^{-2})
\sum_{k,l}(-q^2)^{l-k}x_{\mu,k}A_{lk}x_{l,\nu}.\] We use induction
on $m$ to prove it. If $m=1$, we have
\[x_{11}x_{i+1,j+1}=x_{i+1,j+1}x_{11}+(q^2-q^{-2})x_{i+1,1}x_{1,j+1}\]
which is one of the defining relations. In general, by quantum
Laplace expansion,
\begin{eqnarray*}\label{fterm}
det_qA  x_{\mu, \nu}
&=& \sum_s(-q^2)^{m-s}A_{ms}x_{ms}x_{\mu, \nu}\\
&=&
\sum_s(-q^2)^{m-s}A_{ms} \left [x_{\mu, \nu}
x_{ms}+(q^2-q^{-2})x_{\mu,s}x_{m,\nu} \right ]\\
&=&\sum_s(-q^2)^{m-s}A_{ms}x_{\mu, \nu}x_{ms}\\
&&+ (q^2-q^{-2})\sum_s (-q^2)^{s-m}x_{\mu,s}A_{ms}x_{m, \nu}.
\end{eqnarray*}
Denote by $A_{m,l;s,k}$ the determinant of the quantum matrix
obtained from the quantum matrix $A$ by deleting the $m, l$-th
rows and $s, k$-th columns. It follows from the induction
hypothesis that
\begin{eqnarray*}
&&\sum_s(-q^2)^{m-s}A_{ms}x_{\mu, \nu}x_{ms}\\
&=&\sum_s(-q^2)^{m-s}\left [A_{ms}x_{\mu,
\nu}x_{ms}+(q^2-q^{-2})\sum_{k,l}
 x_{\mu,k}A_{m,l;s,k}x_{l,\nu}x_{ms}\right]
\\
&=&x_{\mu, \nu}det_qA+(q^2-q^{-2})\sum_{k,l}
x_{\mu,k}A_{lk}x_{l,\nu}.\end{eqnarray*} Combining this with the
last equation, we get the desired formulae.

Note that $\bar{D}$ is the inverse matrix of $D^\prime$. Hence,
$det_qA$ commutes with all the entries of $\bar{D}$. By the
relation given in \cite{zhangrb}, we can also see that
$det_{q^{-1}}\bar{D}$ commutes with all of the entries in $A$, and
therefore, $det_{q^{-1}}D^\prime =(det_{q^{-1}}\bar{D})^{-1}$
commutes with all of the entries in $A$. The other formulas can be
proved similarly. Consequently,
$det_qA(det_{q^{-1}}D^\prime)^{-1}$ is a central element.

\qed

The following result shows that the constructions of the quantized
function algebras of the quantum general supergroup in \cite{f},
\cite{manin1} and \cite{zhangrb} are in fact equivalent.

\begin{Thm}\label{equivalent}
As Hopf algebras,
\[G_q\cong {\s O}_q(GL_{m\mid n}).\]
\end{Thm}

\pof Identifying $T^{(\lambda)}$ with $L(\lambda)\otimes
L(\lambda^\dagger)$ as $U_q(\mathfrak g)$-bimodules, we obtain
\[G_q^{\pi}\cong \bigoplus_{\lambda\in \Lambda^{(1)}}
L(\lambda)\otimes L(\lambda^\dagger), \] upon using quantum
Peter-Weyl theorem, Proposition~\ref{peter-weyl}.

By the universal property of ${\s O}_q(M_{m\mid n})$
\cite{manin1}~Theorem 1.6, there is a surjective homomorphism
\begin{eqnarray*}\Phi: {\s O}_q(M_{m\mid n})\longrightarrow G_q^{\pi}, &&
x_{ij} \mapsto
 (-1)^{[i][j]+[j]}q^{(\epsilon_i,\epsilon_i)-(\epsilon_j,\epsilon_j)}
t_{ij}, \\
&& i,j=1,2,\cdots, m+n.\end{eqnarray*} The map also preserves the
co-product and co-unit, as can be easily seen by inspection.

The algebra ${\s O}_q(M_{m\mid n})$ is ${\mathbb Z}_+$-graded with
gradation assignment $degx_{ij}=1$. Denote by ${\s O}_q(M_{m\mid
n})_k$ the homogeneous component of degree $k$. Also note that
$G_q^{\pi}=\bigoplus_k (G_q^{\pi})_k,$ with
\[(G_q^{\pi})_k\cong\bigoplus_{\lambda:|\lambda|=k}L(\lambda)\otimes L(\lambda^\dagger),\]
where $|\lambda|=\sum_a\lambda_a$ for
$\lambda=(\lambda_1,\lambda_2,\cdots,\lambda_m|\lambda_{m+1},\cdots,
\lambda_{m+n})\in \Lambda^{(1)}$. By Proposition 3 of
\cite{zhangjmp}, each $L(\lambda)$ has the same dimension as its
classical counter part, thus by Proposition 3.3 of \cite{sz},
\[dim(G_q^{\pi})_k=\sum_{r=0}^k
\begin{pmatrix}m^2+n^2+r-1\\r
\end{pmatrix}\begin{pmatrix}2mn\\k-r\end{pmatrix}\]
which equals to $dim {\s O}_q(M_{m\mid n})_k$. Hence, as
bialgebras,
\[G_q^{\pi}\cong {\s O}_q(M_{m\mid n}).\]
The natural embedding of $G_q^{\pi}$ in $G_q$ leads to an
embedding of ${\s O}_q(M_{m\mid n})$ in $G_q$.

As an intermediate step to the proof of the theorem, we introduce
the localization of ${\s O}_q(M_{m\mid n})[(det_q A)^{-1}]$ of ${\s
O}_q(M_{m\mid n})$ at $det_q A$. Then $A$ is invertible in ${\s
O}_q(M_{m\mid n})[(det_q A)^{-1}]$, thus the entries of the matrix
$D-CA^{-1}B$ all belong to this superalgebra. We shall still
denote this matrix by $D^\prime$ by an abuse of notation. Now we
localize ${\s O}_q(M_{m\mid n})[(det_q A)^{-1}]$ at
$det_{q^{-1}}D^\prime$, and denote the resulting superalgebra by
${\s O}^\prime_q(GL_{m\mid n})$. Obviously ${\s
O}^\prime_q(GL_{m\mid n})$ is isomorphic to ${\s O}_q(GL_{m\mid
n})$.

Now we want to show that $G_q$ and ${\s O}^\prime_q(GL_{m\mid n})$
are isomorphic.

In the superalgebra $G_q$, the quantum supermatrix $(t_{ij})$ is
invertible with the inverse matrix $(\bar{t}_{ij})$. The antipode
maps the quantum supermatrix to its inverse matrix. Hence,
$\bar{t}_{ij}$'s can be obtained from $t_{ij}$'s together with
$(det_q A)^{-1}$ $(det_{q^{-1}}D)^{-1}$. Therefore, we have a
surjection from ${\s O}_q(GL_{m\mid n})$ to $G_q$, which induces a
surjective map $\phi: {\s O}^\prime_q(GL_{m\mid n})\rightarrow
G_q$.

Let $k$ belong to the kernel of the map $\phi$. Then by using
Proposition \ref{commun} and the fact that $det_qA$ and
$det_{q^{-1}}D^\prime$ commute, we see that for some positive
integer $i$ and a sufficiently large $j$, \[ (det_q A)^i
(det_{q^{-1}}D^\prime)^j k \in {\s O}_q(M_{m\mid n}). \] This
element, belonging to $ker\phi$, must vanish. However, the quantum
determinants $det_q A$ and $det_{q^{-1}}D^\prime$ are invertible
in in ${\s O}^\prime_q(GL_{m\mid n})$, thus $k=0$. This proves the
injectivity of $\phi$, thus establishing that ${\s O}_q(GL_{m\mid
n})$ and $G_q$ are isomorphic as associative superalgebras.

The fact that the co-algebraic structures of ${\s O}_q(GL_{m\mid
n})$ and $G_q$ also coincide follows from computations in
\cite{f}, which showed that $\Delta(det_q A)$ and
$\Delta(det_{q^{-1}}D)$ are invertible elements of $G_q\otimes
G_q$.

\qed

\section{Construction of bases}\label{BASES}

We shall follow \cite{lusztig} to construct bases for the Hopf
superalgebras ${\s O}_q(GL_{m\mid n})$ and  ${\s O}_q(SL_{m\mid
n})$. We expect the bases to be useful for studying the structure
of these Hopf superalgebras. In \cite{manin1}, it was proved that
the ordered monomials form a basis for the superalgebra ${\s
O}_q(M_{m\mid n})$. However, no results seem to be available in
the literature on bases for ${\s O}_q(GL_{m\mid n})$ and  ${\s
O}_q(SL_{m\mid n})$.  As we have pointed out in the previous
section,  some quantum minors behave very well under the left and
right actions of the quantum enveloping superalgebra
$U_q(\mathfrak g)$. Hence, it is natural to expect that any nice
basis of ${\s O}_q(GL_{m\mid n})$ should contain these quantum
minors. Also, if a basis of ${\s O}_q(GL_{m\mid n})$ is invariant
under the multiplication of the quantum Berezinian, then we can
get from it a basis for ${\s O}_q(SL_{m\mid n})$ by setting the
quantum Berezinian to $1$.

The map defined in the lemma below is inspired by the definition
of the anti-automorphism $\bar{ }$ \ in the previous section and
hence will be denoted by the same natation. It is a main
ingredient for the construction of the dual canonical bases.
\begin{Lem}
\begin{enumerate}
\item
The mapping
\begin{eqnarray*}^-:x_{ij}\mapsto x_{ij}, \quad
q\mapsto q^{-1}\end{eqnarray*} extends to a superalgebra
anti-automorphism of ${\s O}_q(M_{m\mid n})$ regarded as a
superalgebra over ${\mathbb Q}$.
\item The anti-automorphism $\bar{ }$ of
${\s O}_q(M_{m\mid n})$ extends uniquely to
${\s O}_q(GL_{m\mid n})$ by requiring
\[\overline{(det_qA)^{-1}}=(det_qA)^{-1},
\quad \overline{(det_{q^{-1}}D)^{-1}}=(det_{q^{-1}}D)^{-1}.\]
\end{enumerate}
\end{Lem}
The lemma can be proved easily by inspecting the defining
relations.

\begin{Rem} The anti-automorphism $^-$ commutes with
the isomorphism $\phi$ in the proof of Theorem \ref{equivalent}.
Indeed, one can directly check that the elements
$(-1)^{[i][j]+[j]}q^{(\epsilon_i,\epsilon_i)-(\epsilon_j,\epsilon_j)}
t_{ij}$ are bar invariant by using Proposition~\ref{dagger} and
Lemma~\ref{anti}.\end{Rem}

Arrange the generators according to the lexicographic order, namely
\begin{eqnarray*}& &x_{11}< x_{12}<\cdots x_{1 m+n}<x_{21}<\cdots \\\nonumber
& &<x_{m, 1}<x_{m, 2}<\cdots <x_{m, m}<\cdots<x_{m+n,
m+n}.\end{eqnarray*} For any matrix $M=(m_{ij})\in
M_{m+n}({\mathbb Z}_+)$, $m_{ij}=0,1$ if $[i]+[j]=\bar{1}$, we
define a monomial $x^M$ by
\begin{equation} x^M=x_{11}^{m_{11}}x_{12}^{m_{12}}\cdots x_{1m+n}^{m_{1m+n}}
x_{21}^{m_{21}}\cdots x_{2m+n}^{m_{2m+n}} \cdots
x_{m+n,m+n}^{m_{m+n,m+n}}. \end{equation} 
Observe that the factors are
arranged in the lexicographic order.

To construct a basis using Lusztig's method \cite{lusztig}, we
need to modify the monomials. Define the {\em normalized
monomials}
$$x(M)=q^{-\sum_{i,
j<k}(-1)^{[i]}m_{ij}m_{ik}-\sum_{l,
s<t}(-1)^{[l]}m_{sl}m_{tl}}x^M. $$ We shall impose a partial order
on the set of the normalized monomials by given a partial order to
the matrices $M$ in the following way. Let $M=(m_{ij})\in
M_{m+n}({\mathbb Z}_+)$. If $m_{ij}m_{st}\ge 1$ for two pairs of
indices $i,j$ and $s,t$ satisfying $i<s, j<t$, we define a new
matrix $M^\prime =(m_{uv}^\prime)\in M_{m+n}({\mathbb Z}_+)$ with
\begin{eqnarray*}
m_{ij}^\prime=m_{ij}-1, && m_{st}^\prime=m_{st}-1,\\
m_{it}^\prime=m_{it}+1, && m_{sj}^\prime=m_{sj}+1,\\
m_{uv}^\prime=m_{uv},  && \mbox{ for all other entries}.
\end{eqnarray*}
We say that the matrix $M^\prime$ is obtained from the matrix $M$
by a $2\times 2$ sub-matrix transformation. Using this we may
define a partial order on the set $M_{m+n}({\mathbb Z}_+)$ such
that $M< N$ if $M$ can be obtained from $N$ by a sequence of
$2\times 2$ sub-matrix transformations.

Given $M=(m_{ij})\in M_{m+n}({\mathbb Z}_+)$, we define
the row sums $ro(M)$ and the column sums $co(M)$ of the matrix,
respectively, by
\begin{eqnarray*}
ro(M)&=&(\sum_j m_{1j},\cdots,\sum_jm_{m+n,j})=(r_1(M),r_2(M),\cdots,r_{m+n}(M)), \\
co(M)&=&(\sum_j m_{j1},\cdots,\sum_j
m_{j,m+n})=(c_1(M),c_2(M),\cdots,c_{m+n}(M)).
\end{eqnarray*}
Note that the $2\times 2$ sub-matrix transformations keep the row
sums and column sums unchanged. Let us also introduce the
following notation, which will be frequently used below:
\[
\M=\left\{\left.\begin{pmatrix}M_1 & M_2\\
M_3& M_4 \end{pmatrix}\right| \begin{array}{l l} M_1\in
M_m({\mathbb Z}_+), &M_2\in M_{m\times n}({\mathbb Z}_2), \\
M_3\in M_{n\times m}({\mathbb Z}_2), & M_4\in M_n({\mathbb
Z}_+)\end{array}\right\}.
\]
Whenever an element of $\M$ is considered, we assume that it is in
this block form.

The following lemma will be needed when constructing the dual
canonical basis. It follows directly from the defining
relations of the algebra ${\s O}_q(M_{m\mid n})$.

\begin{Lem}\label{bar}For any $M\in \M$,
\[\overline{x(M)}=x(M)+\sum_{T<M}c_{M,T}x(T),\]
where the coefficients $c_{M,T}\in {\mathbb Z}[q,q^{-1}]$.\end{Lem}

We first construct a basis for the subalgebra $H$ generated by the
entries of the matrices $A$ and $B$, namely, the entries $x_{ij},
i\le m$ of the quantum supermatrix.  We have the following result.
\begin{Thm}\label{AB} For any
$\begin{pmatrix}M_1&M_2\\0&0\end{pmatrix}\in \M$, there exists a
unique element $\Omega_q\begin{pmatrix}M_1&M_2\\0&0\end{pmatrix} $
determined by the following conditions:
\begin{enumerate}
\item $\overline{\Omega_q\begin{pmatrix}M_1&
M_2\\0&0\end{pmatrix}}=\Omega_q\begin{pmatrix}M_1&M_2\\0&0\end{pmatrix}$.
\item \label{expression}
$\Omega_q\begin{pmatrix} M_1&M_2\\0&0\end{pmatrix}
=x\begin{pmatrix}M_1&M_2\\0&0\end{pmatrix}+ \sum
h_{M_1M_1^\prime}^{M_2M_2^\prime}x\begin{pmatrix}M_1^\prime
&M_2^\prime\\0&0\end{pmatrix}$,  where $h_{M_1 M_1^\prime}^{M_2
M_2^\prime}$ $\in$  $q{\mathbb Z}[q]$. Here the sum is over
$\begin{pmatrix}M_1^\prime &M_2^\prime\\0&0\end{pmatrix} \in \M$
satisfying the condition $\begin{pmatrix}M_1^\prime
&M_2^\prime\\0&0\end{pmatrix} <
\begin{pmatrix}M_1&M_2\\0&0\end{pmatrix}$.
\end{enumerate}
The elements $\Omega_q\begin{pmatrix}M_1&M_2\\0&0\end{pmatrix}$
form a basis of $H$.

The quantum minors $\Delta((i_1,i_2,\cdots, i_r),
(j_1,j_2,\cdots,j_r))^*$ with $1\le i_1< i_2<\cdots< i_r$, $1\le
j_1<j_2\cdots<j_r\le m+n$ are basis elements.
\end{Thm}

\pof We only need to prove the last statement. For simplicity, we
shall only consider the statement for $\Delta([1,r], [s,s+r])^*$.
From the definition we can see that $\Delta([1,r], [s,s+r])^*$ are
indeed of the form as that in (\ref{expression}). Now, we show the
bar invariance of these quantum minors. It is clear that
\[\delta^*\circ \psi =(-\otimes \psi)\circ\delta^*,\]
where $\psi$ is the anti-automorphism of $A_q^*$ fixing all
generators and sending $q$ to $q^{-1}$. Note that both sides are
algebra anti-automorphisms so we only need to check the
generators. Another fact we need is
\[\psi(\xi_s\xi_{s+1}\cdots\xi_{s+r})=
q^{r(r+1)}\xi_s \xi_{s+1} \cdots \xi_{s+r}.\]
Hence,
\begin{eqnarray*}\delta^*(\psi(\xi_1 \xi_2 \cdots \xi_r))
&=&\overline{\Delta([1,r], [s,s+r])^*}\otimes
\psi(\xi_s\xi_{s+1}\cdots\xi_{s+r})+\cdots\\\nonumber
&=&q^{r(r+1)}\overline{\Delta([1,r], [s,s+r])^*}\otimes
\xi_s\xi_{s+1}\cdots\xi_{s+r}+\cdots,\end{eqnarray*} which implies
$\overline{\Delta([1,r], [s,s+r])^*}=\Delta([1,r], [s,s+r])^*$.

\qed

For any matrix $M$, denote by $S(M)$ the sum of all entries in
$M$. The quantum determinant $det_qA$ of the quantum matrix $A$ is
a special quantum minor $\Delta([1, m], [1, m])^*$ which has the
following property.
\begin{Lem}\label{A}
\[det_qA \Omega_q\begin{pmatrix}M_1&M_2\\0&0\end{pmatrix}
=q^{S(M_2)}\Omega_q\begin{pmatrix}M_1+I_m&M_2
\\0&0\end{pmatrix}.\]
\end{Lem}

\pof It is known that $det_qA x_{ij}=x_{ij}det_qA$ for
$i,j=1,2,\cdots, m$. The same argument as  Lemma 3.3 in
\cite{jz} shows that  $det_qA x_{i, \mu}=q^2x_{i, \mu}det_qA$, for
$\mu=m+1,m+2,\cdots, m+n$. The relations together with part (1) of Theorem~\ref{AB}
imply that
$q^{-S(M_2)}det_qA\Omega_q\begin{pmatrix}M_1&M_2\\0&0\end{pmatrix}$
is bar invariant. By Theorem 5.2 of \cite{zhanghc},
$q^{-S(M_2)}det_qA\Omega_q
\begin{pmatrix}M_1&M_2\\0&0\end{pmatrix}$ is of the form
$\Omega_q\begin{pmatrix}M_1^\prime &M_2^\prime \\0&0\end{pmatrix}$
which must be equal to
$\Omega_q\begin{pmatrix}M_1+I_m&M_2\\0&0\end{pmatrix}$.

\qed

We can perform similarly analysis for the subalgebra generated by
the entries of $C$ to prove the following result.

\begin{Thm}\label{C} For any $(m+n)\times
(m+n)$ matrix $\begin{pmatrix}0&0\\M_3&0\end{pmatrix}\in\M$, there
exists a unique element
$\Omega_{q^{-1}}\begin{pmatrix}0&0\\M_3&0\end{pmatrix}$ with
properties

\begin{enumerate}
\item $\overline{\Omega_{q^{-1}}\begin{pmatrix}0&0\\M_3&0\end{pmatrix}}
=\Omega_{q^{-1}}\begin{pmatrix}0&0\\M_3&0
\end{pmatrix}$.
\item $\Omega_{q^{-1}}\begin{pmatrix}0&0\\M_3&0\end{pmatrix}
=x\begin{pmatrix}0&0\\M_3&0\end{pmatrix}+\sum_{T_3<M_3}
h_{T_3M_3}x\begin{pmatrix}0&0\\T_3&0\end{pmatrix}$,\\ where
$h_{T_3M_3}\in q^{-1}{\mathbb Z}[q^{-1}]$.
\end{enumerate}
The elements
$\Omega_{q^{-1}}\begin{pmatrix}0&0\\M_3&0\end{pmatrix}$ form a
basis of the subalgebra generated by the entries of $C$. The
quantum minors $\Delta((i_1,i_2,\cdots, i_r), (j_1, j_2,\cdots,
j_r))$ with $m+1\le i_1<i_2<\cdots<i_r\le m+n, 1\le
j_1<j_2<\cdots<j_r\le n$
 are basis
elements for $r\le min\{m,n\}$.\end{Thm}

Now we consider the subalgebra generated by the entries of $A, B,
C$. It has a basis
\[\{N\begin{pmatrix}M_1&M_2\\M_3&0\end{pmatrix}
:=q^{-\sum_ic_i(M_1)c_i(M_3)}\Omega_q\begin{pmatrix}M_1&M_2\\0&0
\end{pmatrix}
\Omega_{q^{-1}}\begin{pmatrix}0&0\\M_3&0\end{pmatrix}\},\] where
$c_i(M_1)$ and $c_i(M_3)$ are the $i$th column sums of $M_1$ and
$M_3$ respectively. The basis is ordered according to the order of
the matrices $\begin{pmatrix}M_1&M_2\\M_3&0\end{pmatrix}$. From
the construction, it is clear that
\[N\begin{pmatrix}M_1&M_2\\M_3&0\end{pmatrix}
:=x\begin{pmatrix}M_1&M_2\\M_3&0\end{pmatrix}+\text{ lower terms}.\]
Hence, by Lemma~\ref{bar}, we have
\begin{Lem}
\[\overline{N\begin{pmatrix}M_1&M_2\\M_3&0\end{pmatrix}}
=N\begin{pmatrix}M_1&M_2\\M_3&0\end{pmatrix}+ \text{ lower
terms}.\]
\end{Lem}
This leads to the following result.
\begin{Thm}
For any given $\begin{pmatrix}M_1&M_2\\M_3&0\end{pmatrix}\in \M$,
there exists a unique element
$\Omega_q\begin{pmatrix}M_1&M_2\\M_3&0\end{pmatrix} $ determined
by the following conditions:
\begin{enumerate}
\item $\overline{\Omega_q\begin{pmatrix}M_1&M_2\\M_3&0\end{pmatrix}}
=\Omega_q\begin{pmatrix}M_1&M_2\\M_3&0\end{pmatrix}$.
\item
$\Omega_q\begin{pmatrix}M_1&M_2\\M_3&0\end{pmatrix}
=N\begin{pmatrix}M_1&M_2\\M_3&0\end{pmatrix}+ \sum h_{M_1 M_2
M_3}^{M_1^\prime M_2^\prime M_3^\prime} N\begin{pmatrix}M_1^\prime
&M_2^\prime\\M_3^\prime&0\end{pmatrix}$ where the summation is
over all the matrices $\begin{pmatrix}M_1^\prime
&M_2^\prime\\M_3^\prime&0\end{pmatrix}<
\begin{pmatrix}M_1&M_2\\M_3&0\end{pmatrix}$,
and $h_{M_1 M_2 M_3}^{M_1^\prime M_2^\prime M_3^\prime} \in
q{\mathbb Z}$$[q]$.
\end{enumerate}
The elements $\Omega_q\begin{pmatrix}M_1&M_2\\M_3&0\end{pmatrix}$
form a basis of the subalgebra generated by the entries of $A, B,
C$.
\end{Thm}

Arguments analogous to Theorem 4.3 in \cite{jz} show that $det_qA$
$q$-commutes with all entries in $C$. Hence, in a way similar to
the proof of Lemma~\ref{A}, we can show that
\begin{Lem}\label{AA} {\ \ }
$det_qA \Omega_q\begin{pmatrix}M_1&M_2\\M_3&0\end{pmatrix}
=q^{S(M_2)+S(M_3)}\Omega_q
\begin{pmatrix}M_1+I_m&M_2\\M_3&0\end{pmatrix}.$
\end{Lem}

Recall the definition of the matrix $D'$. We have
\begin{Prop}The entries of  the matrix $D^\prime=(y_{\mu, \nu})
:=D-CA^{-1}B$ are bar invariant.\end{Prop}
\pof The entries of $D^\prime$ are of the form:
\[x_{\mu, \nu}-\sum_{k,l=1}^m(-q^2)^{l-k}
x_{\mu, k}A_{lk}(det_qA)^{-1}x_{l,\nu},\] where $A_{lk}$ is the
quantum minor obtained from $A$ by deleting the $l$th row and
$k$th column which is bar invariant by \cite{zhanghc}~Lemma 3.3.

Since $x_{\mu, \nu}$ and $(det_qA)^{-1}$ are bar invariant and
$(det_qA)^{-1}$ $q$-commutes with $x_{\mu, k}$ and $x_{l,\nu}$, we
only need to show that
\[d_{ij}:=\sum_{k,l=1}^m(-q^2)^{l-k}x_{\mu, k}A_{lk}x_{l,\nu}\]
is bar invariant. By repeatedly using quantum Laplace expansion,
we have
\begin{eqnarray*}\overline{d_{ij}}&=&-\sum_{k,l=1}^m(-q^2)^{k-l}
x_{l,\nu}A_{lk}x_{\mu, k}\\\nonumber
&=&-\sum_{k,l=1}^m(-q^2)^{m-l}(-q^2)^{m-k}x_{l,\nu}x_{\mu,
k}A_{lk}\\\nonumber
&=&\sum_{k,l=1}^m(-q^2)^{m-k}(-q^2)^{m-l}x_{\mu,
k}x_{l,\nu}A_{lk}\\\nonumber
&=&\sum_{k,l=1}^m(-q^2)^{m-k}(-q^2)^{l-m}x_{\mu,
k}A_{lk}x_{l,\nu}\\ &=&d_{ij}.\end{eqnarray*} \qed

For any matrix $M=\begin{pmatrix}0&0\\0&M_4\end{pmatrix}\in \M$,
we define the ordered monomials $y(M)$ in the same way as $x(M)$.
The monomials $y(M)$ form a basis of the subalgebra generated by
the entries of the quantum $q^{-1}$-matrix $D^\prime$. Using the
same method as that in \cite{zhanghc}, we get a basis of the
subalgebra generated by the entries of $D^\prime$.
\begin{Thm}\label{D}
For any matrix $\begin{pmatrix}0&0\\0&M_4\end{pmatrix}\in \M$,
there exists a unique element
$\Omega_{q^{-1}}\begin{pmatrix}0&0\\0&M_4\end{pmatrix}$ with the
properties:
\begin{enumerate}
\item $\overline{\Omega_{q^{-1}}\begin{pmatrix}0&0\\0&M_4\end{pmatrix}}
=\Omega_{q^{-1}}\begin{pmatrix}0&0\\0&M_4\end{pmatrix},$

\item $\Omega_{q^{-1}}\begin{pmatrix}0&0\\0&M_4\end{pmatrix}
=y\begin{pmatrix}0&0\\0&M_4\end{pmatrix}+ \sum_{T_4<M_4}h_{T_4M_4}
y\begin{pmatrix}0&0\\0&T_4\end{pmatrix}$ where the coefficients
$h_{T_4M4}\in q^{-1}{\mathbb Z}[q^{-1}]$.\end{enumerate} the
elements $\Omega_{q^{-1}}\begin{pmatrix}0&0\\0&M_4\end{pmatrix}$
form a basis of the subalgebra generated by the entries of
$D^\prime$. In particular, all the quantum minors of the
$q^{-1}$-matrix $D^\prime$ are basis elements.
\end{Thm}

By Proposition 3.6 in \cite{zhanghc}, we have
\begin{Prop}\label{DD}
\[det_{q^{-1}}D^\prime \Omega_{q^{-1}}
\begin{pmatrix}0&0\\0&M_4\end{pmatrix}
=\Omega_{q^{-1}}\begin{pmatrix}0&0
\\0&M_4+I_n\end{pmatrix}.\]
\end{Prop}

Now we proceed to the construction of a basis of ${\s
O}_q(GL_{m\mid n})$. For any $a, d \in {\mathbb Z}$, and
$M=\begin{pmatrix}M_1&M_2\\M_3&M_4\end{pmatrix}\in\M$,  let
\[\Psi(M; a, d)=\sum_j c_j(M_2)c_j(M_4)-\sum_j
r_j(M_3)r_j(M_4)-(a+d)(S(M_2)+S(M_3)), \] and define
\begin{eqnarray*} N_{a, d}(M)
:&=&q^{\Psi(M; a, d)}(det_qA)^a \Omega_q
\begin{pmatrix}M_1&M_2\\M_3&0\end{pmatrix}
\Omega_{q^{-1}}\begin{pmatrix}0&0\\0&M_4\end{pmatrix}(det_q
{D^\prime})^d.\end{eqnarray*}

\begin{Prop} The elements
\[N_{a,d}\begin{pmatrix}M_1&M_2\\M_3&M_4\end{pmatrix}\]
form a basis of algebra ${\s O}_q(GL_{m\mid n})$, where $a, d\in
{\mathbb Z}$, and
$\begin{pmatrix}M_1&M_2\\M_3&M_4\end{pmatrix}\in\M$  satisfies the
condition that $M_1$ and $M_4$ must have at least one zero
diagonal entry each.
\end{Prop}
\pof
It was proved in \cite{manin1} that the ordered monomials
form a basis of the algebra ${\s O}_q(M_{m\mid n})$. Thus the
following set of elements
\[P(M; a, d):=q^{\Psi(M; a, d)}(det_qA)^a x
\begin{pmatrix}M_1&M_2\\M_3&0\end{pmatrix}
x\begin{pmatrix}0&0\\0&M_4\end{pmatrix}(det_q {D^\prime})^d\]
form a basis of ${\s O}_q(GL_{m\mid n})$, where $a, d\in {\mathbb
Z}$, and $M=\begin{pmatrix}M_1&M_2\\M_3&M_4\end{pmatrix}$ $ \in
\M$ satisfies the condition that $M_1$ and $M_4$ must have at
least one zero diagonal entry each. The order of the monomials
$x(M)$ induces an order on the above basis, where $P(M; a, d)\ge P(M';
a', d')$ if and only if $a<a^\prime$ or $a=a^\prime$ and
$d<d^\prime$ or $a=a^\prime, d=d^\prime$ but $M\ge M^\prime.$ The
element $N_{a, d}\begin{pmatrix}M_1&M_2\\M_3&M_4\end{pmatrix}$ can
be written as
\begin{eqnarray*}N_{a, d}\begin{pmatrix}M_1&M_2\\M_3&M_4\end{pmatrix}
=P(M; a, d)+\text{lower terms}\end{eqnarray*} Hence, the statement
follows.

\qed

For the construction of a basis, we shall need the following lemma
which is derived directly from the defining relations of the
algebra  ${\s O}_q(M_{m\mid n})$.

\begin{Lem}For any $M\in \M$ and $a,d\in{\mathbb Z}$,
\[\overline{N_{a,d}(M)}=N_{a,d}(M)+\sum_{N_{a^\prime, d^\prime}(T)<N_{a,d}(M)}
c_{a,d,a^\prime,d^\prime, M,T}N_{a^\prime,d^\prime}(T)\]
where $c_{a,d,a^\prime,d^\prime, M,T}\in {\mathbb Z}[q, q^{-1}]$.
\end{Lem}

By using the lemma, we can prove the following theorem which is one
of the main results of this paper.
\begin{Thm}\label{basis}
There is a unique basis $B_q^*$ of ${\s O}_q(GL_{m\mid n})$
consisting of elements $\Omega_{q, a,d}(M)$ with  $M=(m_{ij})\in
M_{m+n}({\mathbb Z}_+)$ such that $m_{ij}\in {\mathbb Z}_2$, if
$[i]+[j]=\bar{1},$ and $m_{ii}=m_{\mu\mu}=0$  for some $i\le m,
\mu\ge m+1, a,d\in {\mathbb Z}$, which is determined by the
following conditions:
\begin{enumerate}
\item $\overline{\Omega_{q, a,d}(M)}=\Omega_{q,a,d}(M)$ for all $M$.
\item $\Omega_{q,a,d}(M)=N_{a,d}(M)+\sum_{N_{a^\prime,d^\prime}(T)
<N_{a,d}(M)} h_{a,a^\prime,d,d^\prime}(T,M)
N_{a^\prime,d^\prime}(T),$ \\ where
$h_{a,a^\prime,d,d^\prime}(T,M)\in q{\mathbb Z}$$[q]$.
\end{enumerate}
Similarly, there is a unique basis
\[B_{q^{-1}}^*=\{\Omega_{q^{-1},a,d}(M) | M=(m_{ij})\in
M_{m+n}({\mathbb Z}_+)\}\] with $m_{ij}\in {\mathbb Z}_2 \text{ if
}[i]+[j]=\bar{1}, m_{ii}=m_{\mu\mu}=0 \text{ for some }i\le m,
\mu\ge m+1, a,d\in {\mathbb Z}\}$ of ${\s O}_q(GL_{m\mid n})$
determined by the following conditions:
\begin{enumerate}
\item $\overline{\Omega_{q^{-1},a,d}(M)}
=\Omega_{q^{-1},a,d}(M)$ for all $M$.
\item $\Omega_{q^{-1},a,d}(M)=N_{a,d}(M)+ \sum_{N_{a^\prime,d^\prime}(T)
<N_{a,d}(M)} h_{a,a^\prime,d,d^\prime}(T,M)
N_{a^\prime,d^\prime}(T)$,\\  where
$h_{a,a^\prime,d,d^\prime}(T,M) \in q^{-1}{\mathbb Z}$$[q^{-1}]$.
\end{enumerate}
\end{Thm}

\medskip

We shall refer to both $B_q^*$ and $B_{q^{-1}}^*$ as dual
canonical bases of ${\s O}_q(GL_{m\mid n})$. These bases contain
the quantum minors in Theorems~\ref{AB}, \ref{C} and \ref{D} by
construction. Furthermore, we have the following result.
\begin{Thm}
\label{invariance}The bases $B_q^*$ and $B_{q^{-1}}^*$ are invariant
under the multiplication of the quantum Berezinian.
\end{Thm}

\pof Actually, we can show that
\[N_{a,d}(M)Ber_q=N_{a+1,d-1}(M).\]
Write $M=\begin{pmatrix}M_1&M_2\\M_3&M_4\end{pmatrix}$. Since
$Ber_q$ is central, we have
\begin{eqnarray*}
N_{a,d}\begin{pmatrix}M_1&M_2\\M_3&M_4\end{pmatrix}Ber_q&=&
q^{\Psi(M; a, d)} (det_qA)^a\Omega_q
\begin{pmatrix}M_1&M_2\\M_3&0\end{pmatrix}\\
&&\times \Omega_{q^{-1}}\begin{pmatrix}0&0\\0&M_4\end{pmatrix}
(det_{q^{-1}}D^\prime)^d Ber_q\\
&=&q^{\Psi(M; a, d)} (det_qA)^a\Omega_q
\begin{pmatrix}M_1&M_2\\M_3&0\end{pmatrix}\\
&&\times det_qA (det_{q^{-1}}D^\prime)^{-1}
\Omega_{q^{-1}}\begin{pmatrix}0&0\\0&M_4\end{pmatrix}
(det_{q^{-1}}D^\prime)^{d}.\end{eqnarray*} Now, the theorem
follows from Lemma~\ref{AA} and Proposition~\ref{DD} and an
elementary computation of the power of $q$.

\qed

Setting the quantum Berezinian to $1$, we get the superalgebra
${\s O}_q(SL_{m\mid n})$ of functions on the quantum special
linear supergroup $SL_{m\mid n}$, i.e.
\[{\s O}_q(SL_{m\mid n})={\s O}_q(GL_{m\mid n})/<Ber_q-1>,\]
where $<Ber_q-1>$ is the ideal of ${\s O}_q(GL_{m\mid n})$
generated by the central element $Ber_q-1$. Denote by
$M_{m+n}({\mathbb Z}_+)^0$ the subset of $M_{m+n}({\mathbb Z}_+)$
consists of the matrices $M=(m_{ij})_{i,j=1}^{m+n}$ such that
$m_{ss}=m_{rr}=0$ for some $1\le s\le m, m+1\le r\le m+n$ and
$m_{\mu, j}, m_{i,\nu}=0,1$ for all $i,j\le m$, $\mu, \nu\ge m+1$.
Clearly, we get a basis of ${\s O}_q(SL_{m\mid n})$ indexed by
\[{\mathbb Z}\begin{pmatrix}I_m&0\\0&0\end{pmatrix}+M_{m+n}
({\mathbb Z}_+)^0\cup M_{m+n}({\mathbb Z}_+)^0+{\mathbb Z}
\begin{pmatrix}0&0\\0&I_n\end{pmatrix}.\]
The resulting bases are called dual canonical bases of ${\s
O}_q(SL_{m\mid n})$.

We have proved that $\Delta([1,r], [m+n-r+1, m+n])^*$  and
$\Delta([m+n-s+1, m+n], [1, s])$ are dual canonical basis elements
for $r, s\le min\{m,n\}$. We call these quantum minors covariant
quantum minors. The same argument as in \cite{jz}~Theorem 4.3
shows that these covariant quantum minors $q$-commute with all of
the generators $x_{ij}$. Furthermore, similar to the proof as in
\cite{zhanghc}~Theorem 5.2, we can show that
\begin{Thm}\label{covariant} The dual canonical basis $B_q^*$ is
``invariant'' under the multiplication of the covariant minors in
the following sense.
\begin{enumerate}
\item For any dual canonical basis element
$\Omega_{q,a,d}(M)$ (resp. $\Omega_{q^{-1},a,d}(M)$) corresponding
to a matrix $M$ such that the $(i, m+n-i)$ entries are zero for
all $i=1,2,\cdots, r$,
\begin{eqnarray*}
&&\Omega_{q,a,d}(M)\Delta([1, r], [m+n-r+1, m+n])^* \\
&&(\text{resp. }  \Omega_{q^{-1},a,d}(M)\Delta([1, r], [m+n-r+1,
m+n])^*)
\end{eqnarray*}
is also a dual canonical basis element up to a power of $q$.
\item For any dual canonical basis element
$\Omega_{q,a,d}(M)$ (resp. $\Omega_{q^{-1},a,d}(M)$) corresponding
to a matrix $M$ such that the $(m+n-j, j)$ entries are all zero
for all $j=1,2,\cdots, s$,
\begin{eqnarray*}
&&\Omega_{q,a,d}(M)\Delta([m+n-s+1, m+n], [1, s])\\
&& (\text{resp. }  \Omega_{q^{-1},a,d}(M)\Delta([m+n-s+1, m+n],
[1, s])
\end{eqnarray*}
is a dual canonical basis element up to a power of $q$.
\end{enumerate}
\end{Thm}

\section{Invariant subalgebras}\label{INV}

There are two natural  actions of $U_q({\mathfrak g})$ on the
quantized function algebra ${\s O}_q(GL_{m\mid n})$, which
correspond to left and right translations in the classical
setting. However, for convenience, we shall use the algebra anti-automorphism
$\omega: U_q(\mathfrak g)\longrightarrow U_q(\mathfrak g)$ given by
\[\omega(K_i)=K_i,\quad \omega(E_i)=F_i,\quad \omega(F_i)=E_i,\]
to twist the two actions. These two actions are respectively
defined, for all $x\in U_q({\mathfrak g}), f\in {\s O}_q(GL_{m\mid
n})$, by \begin{eqnarray*} R_x(f)&=&\sum_{(f)}f_{(1)}<f_{(2)},
\omega(x)>,\\
L_x(f)&=&\sum_{(f)}<f_{(1)}, \omega(x)>f_{(2)}(-1)^{[x][f_{(1)}]},
\end{eqnarray*}
where we have
used Sweedler's notation $\Delta(f)=\sum f_{(1)}\otimes f_{(2)}$. Note that
$L$ is a left action while $R$ is a right action.
Furthermore, the two  actions commute.

For simplicity, we shall use $x.f$ and $f.x$
to denote the left and right actions respectively.
The actions can be written down explicitly  in terms of generators
as follows:
\[E_i.x_{kl}=\delta_{i,k-1}\delta_{il},\quad F_i.x_{kl}
=\delta{ik}x_{i+1, l},\quad K_i.x_{kl}=q^{2\delta_{ik}}x_{kl}.\]
\[x_{kl}.E_i=\delta_{i+1,l}x_{ki},\quad x_{kl}.F_i
=\delta_{li}x_{k,i-1}, \quad x_{kl}.K_i=q^{2\delta_{li}}x_{kl}.\]

Under these actions, the entries of the matrix $D^\prime$ have the
following property.
\begin{Lem}
For any $\mu,\nu=m+1, m+2,\cdots, m+n$ and $i=1,2,\cdots, m+n-1$,
\[E_i.y_{\mu +1, \nu}=\delta_{i\mu}y_{\mu ,\nu},
\quad F_i.y_{\mu ,\nu}=\delta_{i\mu}y_{\mu +1, \nu},\]
\[y_{\mu , \nu}.F_j=\delta_{j,\nu+1}y_{\mu ,\nu-1},
\quad y_{\mu ,\nu}.E_j=\delta_{j\nu}y_{\mu ,\nu+1}.\]
\end{Lem}
\pof For any $\mu,\nu=m+1,m+2,\cdots,m+n$,
\[y_{\mu , \nu}=x_{\mu ,\nu}-\sum_{k,l=1}^m x_{\mu ,k}
(-q^2)^{l-k}A_{lk}(det_qA)^{-1}x_{l,\nu}.\]
Using the formula for
the coproduct,
\begin{eqnarray*}
\Delta(E_i) &=& E_i \otimes
K_iK_{i +1}^{-1} + 1 \otimes E_i, \\
\Delta(F_i) &=& F_i\otimes 1 + K_i^{-1} K_{i +1}
\otimes F_i,\end{eqnarray*} we have
\begin{eqnarray*}E_i.(x_{\mu +1,\nu}-\sum_{k,l}
x_{\mu +1,k}(-q^2)^{l-k}A_{lk}(det_qA)^{-1}x_{l,\nu})\\\nonumber
=\delta_{i\mu}(x_{\mu ,\nu}-\sum_{k,l}x_{\mu ,k}(-q^2)^{l-k}A_{lk}
(det_qA)^{-1}x_{l,\nu}).\end{eqnarray*}
\begin{eqnarray*}F_i.(x_{\mu ,\nu}-\sum_{k,l}
x_{\mu ,k}(-q^2)^{l-k}A_{lk}(det_qA)^{-1}x_{l,\nu})\\\nonumber
=\delta_{i\mu}(x_{\mu +1,\nu}-\sum_{k,l}x_{\mu +1,k}(-q^2)^{l-k}A_{lk}
(det_qA)^{-1}x_{l,\nu}).\end{eqnarray*} The other formulae can be
proved similarly.

\qed

\begin{Lem}For all $i$,
\[E_i.det_{q^{-1}}D^\prime=det_{q^{-1}} D^\prime .F_i=0, \quad
E_i.Ber_q=Ber_q.F_i=0.\]
\end{Lem}

\pof It is easy to check that $E_m.y_{m+1,\nu}=0$ for all
$\nu=m+1,m+2,\cdots, m+n$. Indeed,
\begin{eqnarray*}E_m.y_{m+1,\nu}&=&x_{m,\nu}-\sum_{k,l=1}^mx_{m,k}
(-q^2)^{l-k}A_{lk}(det_qA)^{-1}x_{l,\nu})\\\nonumber
&=&x_{m,\nu}-\delta_{ml}x_{l,\nu}=0,
\end{eqnarray*}
which implies that $E_m.det_{q^{-1}}D^\prime=0$.

Similarly, we can show that
\[y_{\mu , m+1}.F_{m+1,m}=0.\]
Clearly, $E_i.det_{q^{-1}}D^\prime=0$ for $i\le m$. If $i\ge m$,
$E_i.det_{q^{-1}}D^\prime=0$ is due to the quantum Laplace
expansion. It is known that $E_i.det_qA=0$ for all $i$. This
together with the formula
$Ber_q=det_qA(det_{q^{-1}}D^{\prime})^{-1}$ imply that $E_i.Ber_q=0$ for all $i$.

Similarly, we can show that
$det_qA.F_i=0, det_{q^{-1}}D^\prime .F_i=0, Ber_q .F_i=0$, for all $i$

\qed

\medskip

We shall employ the dual canonical bases constructed to study invariant
subalgebras of ${\s O}_q(GL_{m|n})$ and ${\s O}_q(SL_{m|n})$ under
left and right translations.  Any subset $S$ of the generators
$\{E_i, \, F_i, \, K_i^{\pm 1} \mid i=1, 2, \dots,m+n-1\}$
generates a subalgebra $U_S$ of $U_q({\mathfrak g})$.
\begin{Def} $^{L_S}{\s O}_q(GL_{m\mid n}):=\{f\in
{\s O}_q(SL(n)) \mid x.f=\epsilon(x) f,  \forall x\in
U_S\}$.\end{Def}
It can be easily shown that this is a subalgebra of ${\s
O}_q(GL_{m\mid n})$. It consists of the elements which are
invariant under the left action $L$ of $U_S$. As left and right
translations commute, $^{L_S}{\s O}_q(GL_{m\mid n})$ forms a right
$U_q({\mathfrak g})$-module under $R$. Thus if $T$ is another
subset of the Chevalley generators and the $K_i^{\pm 1}$, we can
also consider \[ {^{L_S}{\s O}_q(GL_{m\mid n})^{R_T}}:=\{f\in
{^{L_S}{\s O}_q(GL_{m|n})} \mid f.x=\epsilon(x) f,  \forall x\in
U_T\}. \] Needless to say, this is a subalgebra of $^{L_S}{\s
O}_q(GL_{m\mid n})$. Below we shall consider in some detail the
subalgebras $n^+$ and $n^-$ (recall that $n^+$ (resp. $n^-$) is
generated by all the $E_i$'s (resp. $F_i$'s)).

Denote by $U_q({\mathfrak g})_0$ the subalgebra of $U_q({\mathfrak
g})$ generated by all even Chevalley generators and the $K_i$'s. Then
\[U_q({\mathfrak g})_0\cong U_q({\mathfrak{gl}}_m)\otimes U_q({\mathfrak{gl}}_n).\]
Let $U_q({\mathfrak p})$ be the subalgebra generated by all
elements of $U_q(\mathfrak g)_0$ and $E_{m}$. For any integral
dominant weight $\lambda$, denote by $L_\lambda^{(0)}$ the
irreducible right $U_q({\mathfrak g})_0$ module with highest
weight $\lambda$. $L_\lambda^{(0)}$ can be extended to a right
$U_q(\mathfrak p)$-module by requiring $E_m$ to act trivially. The
Kac module is the induced module
\[K(\lambda)=Ind_{U_q({\mathfrak p})}^{U_q({\mathfrak g})}
L_\lambda^{(0)}.\] Since $D^\prime$ is a $q^{-1}$-matrix, we can
talk about its quantum minors. In particular, we use
$det_{q^{-1}}D^\prime_s$ to denote the quantum minor of the
$s\times s$ principal sub-matrix of $D^\prime$, for $s\le n$. We
have the following observation.
\begin{Prop} \label{Kac-module}
The subalgebra of invariants ${ }^{L_{n^+}}{\s O}_q (GL_{m\mid
n})^{R_{n^-}}$ is generated by $Ber_q^{\pm 1}$, the quantum minors
$\Delta([1,r], [1, r])^*$ for all
$r=1,2,\cdots, m$,  and $det_{q^{-1}}D^\prime_s$ for all
$s=1,2,\cdots, n$.
\end{Prop}

\pof We can deduce from Theorem 5.2 in \cite{zhangrbj} that
\[{ }^{L_{n^+}}{\s O}_q(GL_{m\mid n})\cong \bigoplus_\lambda K(\lambda),\]
where $\lambda$ ranges over all integral dominant weights. Hence,
the subalgebra ${ }^{L_{n^+}}{\s O}_q(GL_{m\mid n})^{R_{n^-}}$ is
spanned by the lowest weight vectors of all of the Kac modules. It
is known that the lowest weight of each Kac module is of
multiplicity one.

Analogue to the proof of the above lemma, we can see that all
these quantum minors $\Delta([1,r], [1, r])^*$ and  
$det_{q^{-1}}D^\prime_s$ are $L_{n^+}\times R_{n^-}$
invariants. In the previous section it was shown that the
monomials in the quantum minors $\Delta([1,r],
[1, r])^*$, $det_{q^{-1}}D^\prime_s$, and $Ber_q^{\pm
1}$ are linearly independent. Thus in order to prove our claim, we
only need to show that the left weights of these monomials exhaust
all of the integral dominant weights.

Note that the left weights of $\Delta([1,r],
[1, r])^*$, $det_{q^{-1}}D^\prime_s$, and $Ber_q$ are
respectively given by \begin{eqnarray*} (\underbrace{1, \dots,
1}_r, 0, \dots, 0\mid 0, \dots, 0), && (0, \dots, 0\mid
\underbrace{1, \dots, 1}_s, 0, \dots, 0),\\
(1, \dots, 1 \mid -1, \dots, -1).
\end{eqnarray*}
Thus the left weights of their monomials indeed exhaust all the
integral dominant weights.

\qed

\medskip

In the remainder of the paper, we specialize to $n=1$ to study
invariant subalgebras.  To this end, we need to have more detailed
information on $2\times 2$ quantum minors.  For $s\in{\mathbb
Z}_+$, we define
\[[s]_{q^2}=\frac{q^{2s}-1}{q^2-1}, \quad
\begin{pmatrix}s \\ r\end{pmatrix}_{q^2}=\frac{[s]_{q^2}[s-1]_{q^2}
\cdots[s-r+1]_{q^2}}{[r]_{q^2}[r-1]_{q^2} \cdots [1]_{q^2}}.\] For
any indices $1\le j<k\le m+1$, denote
\[M_{jk}:=x_{1j}x_{2k}-q^2x_{1k}x_{2j}.\]
The following lemma will be needed when define Kashiwara operators
which can be proved similarly as the proof of Lemma 2.7 in
\cite{zz}.
\begin{Lem}Assume that $i<k\le m$, $j<l\le m$.
\begin{eqnarray}
&& (x_{ij} x_{kl}-q^2 x_{il}x_{kj})^s=\sum_{m=0}^s (-q^2)^m
\begin{pmatrix}s\\ m\end{pmatrix}_{q^{4}} q^{4m(m-s)}x_{ij}^{s-m}x_{il}^m
x_{kj}^mx_{kl}^{s-m}.\nonumber\\
&&(x_{ij} x_{m+1,m+1}-q^2 x_{i,m+1}x_{m=1,j})^s=\sum_{m=0}^s
(-q^{2})^m
\begin{pmatrix}s\\ m\end{pmatrix}_{q^4} q^{4m(m-s)} \\ \nonumber
&& x_{ij}^{s-m}x_{i,m+1}^m
x_{m+1,j}^mx_{m+1,m+1}^{s-m}.\end{eqnarray}
\end{Lem}

To define the Kashiwara operators $\tilde{E_1}$ and $\tilde{F_1}$,
we need an appropriate basis on which the actions of the Kashiwara
operators are easy to describe.
\begin{Prop} \begin{enumerate}
\item There exists a basis of the algebra ${\s O}_q(GL_{m\mid n})$
consisting  of the elements of the  form
\begin{eqnarray*}&&q^l x\begin{pmatrix}0&\cdots&0&a_{1r}
&a_{1,r+1}&\cdots&a_{1n}\\a_{21}&\cdots&a_{2, r-1}&
a_{2r}&0&\cdots&0\end{pmatrix}\\ \nonumber && \times \Pi M_{ij}
\Pi_{i\ge 3, j} x_{ij}^{a_{ij}},\end{eqnarray*} where
$a_{ij}\in{\mathbb Z}_+$ for all $i,j$, and the $2\times 2$
quantum minors $M_{i j}$ are of the form $\det_q(\{1, 2\}, \ \{i,
j\})$. The product $\Pi M_{ij}$ of quantum minors is arranged
according to the lexicographic order, namely, $M_{ij}\ge M_{st}$
if $j>t$ or $j=t$ and $i\ge s$. 
For given $a_{ij}$'s and the $2\times 2$ minors $M_{ij}$, there is a unique
choice of integer $l$ redering the following property satisfied:
\item The transition matrix between this new basis and the PBW basis
consisting of the modified monomials is of the form:
\begin{equation} \begin{pmatrix}1&\cdots &q{\mathbb Z}[q]\\
0&1&\cdots\\\cdots&\cdots&\cdots\\0&\cdots&1\end{pmatrix}.
\label{matrix-form}
\end{equation}
\end{enumerate}
\end{Prop}

\medskip

Note that when examining the
actions of $\tilde{E}_1$ and $\tilde{F}_1$ on the new basis, we can
ignore the $2\times 2$ minors and those $x_{ij}$ for $i\ge 3$ in
the expression of the new basis elements.
Now the Kashiwara operators $\tilde{E_1}$ and
$\tilde{F_1}$ for the left action  are defined
as follows:
\begin{eqnarray*}&&\tilde{E_1}\left(x\begin{pmatrix}0&0&\cdots&\quad a_{1r} &\cdots&
a_{1n}\\a_{21}&a_{22}&\cdots&a_{2r}&\cdots&0\end{pmatrix}\right) \\
&&=\sum_k q^{\sum_{t=1}^{k-1}2a_{2t}}x\begin{pmatrix}0
&\cdots&\quad 1&\cdots& \quad a_{1r}&\cdots& \quad a_{1n}
\\a_{21}&\cdots&a_{2k}-1&
\cdots&a_{2r}&\cdots&0\end{pmatrix} \end{eqnarray*} where the
summation is over $k$ such that $a_{2k}\ge 1$. Also,
\begin{eqnarray*}
&&\tilde{F_1}\left(x\begin{pmatrix}0&\cdots&a_{1r}&\cdots&a_{1n}\\
a_{21}&\cdots&a_{2r}&\cdots&0\end{pmatrix}\right)\\
&&=\sum_{k}q^{\sum_{t>k}2a_{1t}}x\begin{pmatrix}0&\cdots
&a_{1r}&\cdots&a_{1k}-1&\cdots &a_{1n}\\
a_{21}&\cdots&a_{2r}&\cdots&1&\cdots&0\end{pmatrix},\end{eqnarray*}
where the summation is over $k$ such that $a_{1k}\ge 1$.
Similarly, we can define the Kashiwara operators $\tilde{E_i}, \tilde{F_i}$ for all
$i=1,2,\cdots, m+n-1$.

From the definition of Kashiwara operators and the definition of
$^{L_S}{\s O}_q(G)$, we can show easily that
\begin{Lem} If $E_i, F_j\in S$, then
$$\tilde{E}_i (f) = 0, \quad \tilde{F}_j (f) = 0,
\quad \forall f\in {^{L_S}{\s O}_q(GL_{m\mid n})}, $$ where
$\tilde{E}_i$ and $\tilde{F}_j$ are the Kashiwara operators
associated with $E_i$ and $F_j$.
\end{Lem}

In the following, we let $S$ be any subset of
\[ \{E_i, \, F_i, \, K_i^{\pm 1} \mid i=1, 2, \dots, m+n-1\}\backslash \{F_m\}, \]
and consider the subalgebra of invariants with respect to $S$. We
have the following result.
\begin{Thm}\label{invariant}
The subalgebra of invariants ${^{L_S}{\s O}_q(GL_{m\mid 1})}$ is
spanned by a part of the dual canonical basis $B_q^*$.
\end{Thm}
\pof When $n=1$, the entries of the matrix $C$ $q$-commute with
each other and so the elements
$\Omega_{q^{-1}}\begin{pmatrix}0&0\\M_3&0\end{pmatrix}
=x\begin{pmatrix}0&0\\M_3&0\end{pmatrix}$, for all row vectors
$M_3$ form a basis for the subalgebra generated by entries of $C$.
Furthermore, the basis
\[\left\{N\begin{pmatrix}M_1&M_2\\M_3&0\end{pmatrix}
:=q^{-\sum_ic_i(M_1)c_i(M_3)}\Omega_q
\begin{pmatrix}M_1&M_2\\0&0\end{pmatrix}
\Omega_{q^{-1}}\begin{pmatrix}0&0\\M_3&0\end{pmatrix}\right\},\]
of the subalgebra generated by the entries of the matrices $A$,
$B$ and $C$ is related to the basis
$\left\{x\begin{pmatrix}M_1&M_2\\M_3&0\end{pmatrix}\right\}$ by a
transition matrix of the form (\ref{matrix-form}).  Therefore, the
basis elements
$\Omega_q\begin{pmatrix}M_1&M_2\\M_3&0\end{pmatrix}$ can be
expressed as ${\mathbb Z}[q]$ combinations of the monomials
$x\begin{pmatrix}M'_1&M'_2\\M'_3&0\end{pmatrix}$.

The matrix $D$ has only one element. Thus the basis
$\left\{\Omega_{q, a,
d}\begin{pmatrix}M_1&M_2\\M_3&M_4\end{pmatrix}\right\}$ can be
constructed  as ${\mathbb Z}[q]$ combinations of elements of the
monomial basis
\begin{eqnarray*}q^{-(a+d)(S(M_2)+S(M_3))}(det_qA)^a
x\begin{pmatrix}M_1&M_2\\M_3&M_4\end{pmatrix}
(det_{q^{-1}}D^\prime)^d, \\
\begin{pmatrix}M_1&M_2\\M_3&M_4\end{pmatrix}\in {\mathbb M}, \ a
,d \in {\mathbb Z}.\end{eqnarray*} Now, the same argument as in
\cite{zz} shows that the subalgebra of invariants is spanned by a
part of the dual canonical basis. \qed

Now we consider the invariants with respect to the subalgebra
$U_q(n_+)$ generated by all $E_1, E_2, \cdots, E_{m}$, we deduce
that
\begin{Thm}
Any Kac module is spanned by a part of the dual canonical basis.
\end{Thm}
\pof By Theorem 5.2 in \cite{zhangrbj}, ${ }^{L_{n_+}}{\s
O}_q(GL_{m\mid n})\cong \oplus_\lambda K(\lambda),$ where
$\lambda$ range over all integral dominant weights. Using theorem
\ref{invariant} and consider the left weight space of weight
$\lambda$, we get a  basis of $K(\lambda)$.

\qed

In the case of $GL_{1\mid 1}$, the basis elements are given as:
\[q^{(d-a)(b+c)}x_{11}^a x_{12}^b x_{21}^c (x_{22}+q^2x_{12}x_{11}^{-1}x_{21})^d,\]
where $a,d\in {\mathbb Z}$ and $b,c\in {\mathbb Z}_+$.

In case of $GL_{2\mid 1}$, by Theorem\ref{invariance}, Theorem
\ref{covariant} and the computation for $GL_{1\mid 1}$, we only
need to consider the basis elements parametrized by matrices
\[\begin{pmatrix}a&0&0\\0&0&1\\0&1&0\end{pmatrix},\] where $a\ge 1$
which can be computed easily. Applying the left action of $n^+$,
we get the subalgebra
${ }^{L_{n^+}}{\s O}_q(GL_{2\mid 1})$ which is spanned by
\begin{eqnarray*}
&& q^{-aa_{21}-ba_{12}-a_{11}a_{12}-a_{11}a_{21}-a\alpha-
b\alpha+a\beta+b\beta}(det_q A)^\alpha
x_{11}^{a_{11}}x_{12}^{a_{12}}w^a x_{13}^b
{x_{33}^\prime}^\beta,\\
&& q^{-2\alpha-a_{11}a_{12}-a_{11}-a_{12}+1+2\beta}(det_q
A)^\alpha
x_{11}^{a_{11}}x_{12}^{a_{12}}x_{13}x_{23}{x_{33}^\prime}^\beta,\\
&& q^{-a_{11}a_{12}-a_{12}-2\alpha+2\beta+1}(det_q A)^\alpha
x_{13} w^a{x_{33}^\prime}^\beta\\
&& q^{-a\alpha-\alpha+\beta}x_{11}^{a_{11}}(x_{11}
x_{23}-q^2x_{13}x_{21})^ax_{13}^b
\end{eqnarray*}
where $a,b=0, 1$,  $a_{ij}$ are nonnegative
integers and $\alpha, \beta$ are integers. 
Apply the right actions of $R_{n^-}$, we see that the
subalgebra ${ }^{L_{n^+}}{\s O}_q(GL_{2\mid 1})^{R_{n^-}}$ is
spanned by the following elements: \begin{eqnarray*} &&
q^{\beta-l}x_{12}^l w^a x_{13}{x_{33}^\prime}^\beta, \quad
 q^{b\beta}w^a
x_{13}^b{x_{33}^\prime}^\beta, \quad a, b=0,1,\  l\ge 1, \\
&&q^{2\beta+1-l}x_{12}^lx_{13}x_{23}{x_{33}^\prime}^\beta, \quad
q^{-l}x_{12}^lx_{13}{x_{33}^\prime}^\beta w^a.
\end{eqnarray*}
\begin{Rem} The algebra ${ }^{L_{n^+}}{\s O}_q(GL_{2\mid 1})^{R_{n^-}}$ is
not finitely generated. Indeed, any generating set of the algebra
should contain the elements $x_{12}^l x_{13}$ for all $l\ge 1$.
\end{Rem}

\bibliographystyle{amsplain}

\end{document}